\documentclass[12pt]{amsart}
\usepackage{amscd}
 \pdfoutput=1
 
\usepackage{amsmath}
\usepackage{amsfonts}
\usepackage{amssymb}
\usepackage{amsthm}
\usepackage{mathtools}

\usepackage[utf8]{inputenc}

	\usepackage{pinlabel}
	\usepackage{lmodern}
	\usepackage{textcomp}
	\usepackage[hidelinks]{hyperref}
	\usepackage{color}		
	\usepackage{tabu}
	
\newtheorem{theorem}{Theorem}[section]
\newtheorem{lemma}[theorem]{Lemma}

\newtheorem{definition}[theorem]{Definition}

\theoremstyle{definition}
\newtheorem{remark}[theorem]{Remark}
\newtheorem{example}[theorem]{Example}
\numberwithin{equation}{section}

\newcommand{\Z}{\mathbb{Z}}
\newcommand{\R}{\mathbb{R}}
\newcommand{\C}{\mathbb{C}}

\newcommand{\Ker}{\operatorname{Ker}}
\newcommand{\dbar}{\bar{\partial}}
\newcommand{\M}{\mathcal{M}}
\newcommand{\bm}{\operatorname{bm}}
\newcommand{\inter}{\operatorname{int}}

\newcommand{\triv}{\operatorname{triv}}
\newcommand{\stab}{\operatorname{stab}}
\newcommand{\negg}{\operatorname{ne}}
\newcommand{\pos}{\operatorname{po}}
\newcommand{\tp}{\operatorname{tp}}

\newcommand{\en}{\operatorname{end}}
\newcommand{\U}{\mathcal{U}}
\newcommand{\F}{\mathcal{F}}
\newcommand{\I}{\mathcal{I}}
\newcommand{\Reeb}{\mathcal{R}}
\newcommand{\Pu}{\mathcal{P}}
\newcommand{\Sw}{\mathcal{S}}
\newcommand{\Y}{\mathcal{Y}}

\newcommand{\spn}{\operatorname{Span}}
\newcommand{\Or}{\mathcal{O}}
\newcommand{\capp}{\operatorname{cap}}
\newcommand{\sw}{\text{switch}}
\newcommand{\Flo}{\mathcal{F}}
\newcommand{\pk}{\operatorname{PreKer}}

\begin{document}

\title[To compute orientations of Morse flow trees]
{To compute orientations of Morse flow trees in Legendrian contact homology}
\author{Cecilia Karlsson}

\address{University of Borås,
 Department of Engineering\\
501 90 Borås, Sweden \\
cecilia.karlsson@hb.se}

\begin{abstract}
Let $\Lambda$ be a closed, connected Legendrian submanifold of the 1-jet space of a smooth $n$-dimensional manifold. Associated to $\Lambda$ there is a Legendrian invariant called Legendrian contact homology, which is defined by counting rigid pseudo-holomorphic disks of $\Lambda$. 
Moreover, there exists a bijective correspondence between rigid pseudo-holomorphic disks and rigid Morse flow trees of $\Lambda$, which allows us to compute the Legendrian contact homology of $\Lambda$ via Morse theory. 

If $\Lambda$ is spin, then the moduli space of the rigid disks can be given a coherent orientation, so that the Legendrian contact homology of $\Lambda$ can be defined with coefficients in $\Z$. In this paper we give an algorithm for computing the corresponding orientation of the moduli space of rigid Morse flow trees if the dimension of $\Lambda $ is greater than 1, and up to 4 signs that depend on data that can be extracted from the vertices of the trees, but which are not given explicitly, in the case $n=1$. 
\end{abstract}

\maketitle
\section{Introduction}
Let $M$ be a smooth manifold of dimension $n$, and let $J^1(M)=T^*M \times \R$ denote its 1-jet
space, with  local coordinates $(x_1,\dotsc,x_n,y_1,\dotsc,y_n,z)$. This space can be given the structure  of a contact manifold, with contact distribution $\xi = \Ker ( dz- \sum_i y_i dx_i)$.   
 A submanifold $\Lambda \subset J^1(M)$ is called \emph{Legendrian} if it is $n$-dimensional and everywhere tangent to $\xi$, 
and a \emph{Legendrian isotopy} is a smooth $1$-parameter family of Legendrian submanifolds. The problem of classifying all Legendrian submanifolds of a contact manifold, up to Legendrian isotopy, is a central problem in contact geometry. This motivates finding  Legendrian invariants, that is, invariants that are preserved under Legendrian isotopies. One such invariant is \emph{Legendrian contact homology},  which is a homology theory that fits into the package of Symplectic Field Theory (SFT), introduced by Eliashberg, Givental and Hofer in the paper \cite{sft}. 

More recently, it has been shown that Legendrian contact homology can be used not only to study properties of Legendrian submanifolds, but also to compute symplectic invariants of Weinstein manifolds obtained by surgery along Legendrians, see e.g. ~\cite{surg2,surg1}. It has also been shown, for example in ~\cite{el,rm}, that this can be generalized and used for computations in homological mirror symmetry. This indicates the importance of being able to explicitly understand the Legendrian contact homology complex with integer coefficients. This paper should be an important step in that direction.

Briefly, Legendrian contact homology is the homology of a differential graded algebra (DGA) associated to $\Lambda$.
One should note that Legendrian contact homology has not been worked out in full detail for all contact manifolds, but in the special case when the contact manifold is given by  the 1-jet space $J^1(M) = T^*M \times \R $ of a manifold $M$, then the analytical details were established in \cite{pr}. In the special case $M = \R$, this was also done by Eliashberg in \cite{eliinv} and independently by Chekanov in \cite{chekanov}. 
In the case of a 1-jet space, the DGA of $\Lambda$ can be defined by considering the \emph{Lagrangian projection} 
$\Pi_\C: J^1(M) \to T^*M$. The generators of the algebra are then given by the double points of $\Pi_\C(\Lambda)$, which correspond to \emph{Reeb chords} of $\Lambda$. These are flow segments of the \emph{Reeb vector field}
 $\partial_z$, having start and end point on $\Lambda$.
The differential $\partial$ of the algebra counts certain rigid pseudo-holomorphic disks in $T^*M$ with boundary on $\Pi_\C(\Lambda)$. With a clever choice of almost complex structure on $T^*M$, one gets that the homology of this complex is a Legendrian invariant. 

If $M = \R$, then the count of pseudo-holomorphic disks reduces to combinatorics, as described by Chekanov in \cite{chekanov},  but in higher dimensions the Cauchy-Riemann equations give rise to non-linear partial differential equations, which are hard to solve. 
To simplify a similar problem in Lagrangian Floer homology, Fukaya and Oh in \cite{fo} introduced \emph{gradient flow trees}, which are trees with edges along gradient flow lines in $M$, and gave a one-to-one correspondence between rigid gradient flow trees and rigid pseudo-holomorphic disks in $T^*M$ with boundary on the Lagrangians.
In \cite{trees}, Ekholm generalized this method to also work in the Legendrian contact homology setting.

In \cite{orientbok} it was proved that if $\Lambda$ is spin, then there is a choice of \emph{coherent orientation} of the moduli space of rigid pseudo-holomorphic disks so that Legendrian contact homology can be defined with coefficients in $\Z$. This builds on the work of Fukaya, Oh, Ohta and Ono   in \cite{fooo}, where they describe a way to orient the determinant line of the $\dbar$-operator over the space of trivialized Lagrangian boundary conditions for the unit disk in $\C$. Since the differential in Legendrian contact homology counts rigid pseudo-holomorphic disks, the orientation of the moduli space at a disk $u$ corresponds to a sign of $u$. That the orientation of the moduli space is coherent means that all these signs cancels in $\partial^2$, so that we get $\partial^2 =0$.

In the special case when $\Lambda$ is the conormal lift of a knot $K \subset \R^3$, it was shown in \cite{kch} that there is a coherent orientation scheme of the moduli space of rigid Morse flow trees of $\Lambda$. Moreover, an explicit formula for computing the signs of the trees was derived, and it was proved that this formula in fact coincides with Ng's combinatorially defined Knot contact homology from ~\cite{ ngI, ngII, ngframed}. 
All this has been proved very useful for the aim to study the topology of the underlying knot. In fact, generalizations of the Legendrian contact homology for the conormal lifts of knots gives a complete knot invariant for the underlying smooth knot in $\R^3$, see \cite{enspub}.

 In
the present paper we give an algorithm for orienting the moduli space of rigid Morse flow trees  associated to a closed, connected, spin Legendrian $\Lambda$ in the 1-jet space of a general $n$-dimensional manifold in the case $n>1$. That is, given a rigid Morse flow tree $\Gamma$ of $\Lambda$ we describe a way to compute the orientation of $\Gamma$. In the case when the manifold is $1$-dimensional this algorithm is not completely worked out, but includes 4 implicit signs that depends on combinatorial data of the vertices of $\Gamma$ (and not of $\Gamma$ itself). However, in this case one can as well compute the sign with the results of \cite{sabloffng} or \cite{ngcomp}, for example. 

It should be mentioned that the orientation of $\Gamma$ depends on several choices, which we call \emph{initial orientation choices} and which we describe in Section \ref{sec:triv}. We also point out that the present paper should be considered as a user's guide, and we refer to \cite{orienttrees} for a complete exposition, both for definitions and for proofs.

We summarize the sign formula in the following theorem. 

\begin{theorem}\label{thm:main}
Let $\Lambda$ be a closed, connected, spin Legendrian submanifold of $J^1(M)$, and assume that we have fixed all initial orientation choices. Let $\M$ be the  moduli space of rigid Morse flow trees of $\Lambda$. Then there is a coherent orientation scheme of $\M$ so that the sign $\sigma(\Gamma)$ of a tree $\Gamma \in \M$ is given by
\begin{equation}\label{eq:main}
\sigma(\Gamma)= \nu_{\triv}(\Gamma)\cdot \nu_{\inter}(\Gamma) \cdot  \nu_{\en}(\Gamma)  \cdot \nu_{\stab}(\Gamma),
\end{equation}
where $\nu_{\triv}(\Gamma)$ , $\nu_{\inter}(\Gamma)$ and $\nu_{\en}(\Gamma)$ can be computed completely in terms of $\Gamma$ for all $n$, and the sign $\nu_{\stab}(\Gamma)$ can be computed explicitly in the case $n>1$. 
Moreover, if we use this orientation scheme to define Legendrian contact homology over $\Z$, then we get a Legendrian invariant which is isomorphic to the one defined in \cite{orientbok}.
\end{theorem}

\begin{remark}
The signs $\nu_{\triv}(\Gamma)$ , $\nu_{\inter}(\Gamma)$, $\nu_{\en}(\Gamma)$
and $\nu_{\stab}(\Gamma)$ are given in Definition \ref{def:triv}, \ref{def:nui}, \ref{def:musign} and  \ref{def:stab}, respectively.
In the case $n=1$ the latter sign depends on combinatorial data associated to the vertices of $\Gamma$ (and is not explicitly given as a sign). 
\end{remark}

This theorem follows from the constructions in \cite{orienttrees} together with the results in \cite{teckcob}. Namely, in \cite{orienttrees} the space of rigid pseudo-holomorphic disks of $\Lambda$ is given an orientation which by [\cite{teckcob}, Theorem 2.7] gives rise to a Legendrian contact homology of $\Lambda$ which is invariant under Legendrian isotopies and which is isomorphic to the one in \cite{orientbok}. In \cite{orienttrees} it is then proven that this orientation scheme can be computed in terms of Morse flow trees, via the formula \eqref{eq:main}. In this paper we give a detailed description of this formula.

\begin{remark}
 The orientation scheme in \cite{orientbok} gives the shading rule from \cite{sabloffng} in the case when $M = \R$ and $\Lambda = S^1$.
\end{remark}

 \subsection*{Outline}
To state the definitions that occur in Theorem \ref{thm:main} we first need to introduce some more materials about flow trees, which will be done in Section \ref{sec:trees}. In Section \ref{sec:triv} we give the formula for $\nu_{\triv}$, which is a sign that depends on the trivialization of the tangent bundle of $\Pi_{\C}(\Lambda)$ along the cotangent lift of $\Gamma$, and in particular on the spin structure of $\Lambda$. In Section \ref{sec:intersection} we give the formula for $\nu_{\inter}$, which is a geometric intersection sign coming from intersections of flow manifolds in $M$, and in Section \ref{sec:end} we give the formula for the sign $\nu_{\en}$, which comes from the fact that we treat a special type of vertices as marked points. In Section \ref{sec:stab} we give a formula for $\nu_{\stab}$ in terms of combinatiorial data coming from the tree. 
 This sign is an analytic sign coming from stabilizations of certain Fredholm operators and orientations of dbar problems defined on a disk with trivialized Lagrangian boundary conditions. It is the latter problem which is not completely solved in the case when $n=1$. We finish the paper by giving some examples of sign computations in Section \ref{sec:ex}.

\subsection*{Acknowledgements} 
The author would like to thank Tobias Ekholm and Michael Sullivan for sharing their insights in the subject, and for pointing out mistakes in an earlier version of the paper. She also like to thank the referee for all help with finding mistakes, and for the tremendous work of writing detailed referee reports.

\section{Morse flow trees}\label{sec:trees}
In this section we give a brief introduction to Morse flow trees, following the lines of \cite{trees} and \cite{orienttrees}. We refer the reader to these papers for a more comprehensive description.

Let $\Lambda \subset J^1(M)$ be a closed, spin Legendrian submanifold.
 We will assume that $\Lambda$ is \emph{front generic}, meaning that the base projection $\Pi:J^1(M) \to M$ restricted to $\Lambda$ is an immersion outside a co-dimension 1 submanifold $\Sigma \subset \Lambda$. The points in $\Sigma$ we call \emph{cusp points}. If $\dim(\Lambda)\neq 2 $  we will also assume that $\Lambda$ has \emph{simple front singularities}, meaning that the points in $\Sigma$ are projected to standard cusp singularities under the base projection. If $\dim(\Lambda)=2$ we also allow \emph{swallow-tail singularities}. See  [\cite{trees}, Section 2.2].  
 
 \subsection{Flow lines and vertices}
 Locally the Legendrian $\Lambda$ can be described as the multi-1-jet graph of locally defined functions $f_1,\dotsc,f_l: U \to \R$ for some $U \subset M$ which we will assume to be connected. Each $f_i$ determines a \emph{sheet} of $\Lambda$, given by the 1-jet lift $\{(df_i(x),f_i(x))\subset T^*M\times \R  |x \in U\}$ of $f_i$.  We call such a function $f_i:U \to \R$ a \emph{locally defining function} for $\Lambda$.
 
Fix a metric $g$ on $M$, and let $f_i, f_j$ be locally defining functions for $\Lambda$, defined on $U \subset M$. Then we can consider the function difference $f_i - f_j$ on $U$, which we call a \emph{locally defined function difference associated to $\Lambda$}. Each such function difference gives rise to a negative gradient flow $- \nabla (f_i-f_j)$, where we assume that $f_i \geq f_j$. The solution curves of such a vector field will be called \emph{local flow lines}. If $\gamma(t)$ is such a local flow line we define the \emph{1-jet lift of $\gamma(t)$} to be the pair of curves $\gamma^i(t), \gamma^j(t) \subset J^1(M)$ parametrized by $\gamma^i(t)=(df_i(\gamma(t)),f_i(\gamma(t)))$ and $\gamma^j(t)=(df_j(\gamma(t)),f_j(\gamma(t)))$. The \emph{cotangent lift of $\gamma$} is the pair of curves $\bar \gamma^i, \bar \gamma^j$ given by the Lagrangian projection of $\gamma^i$ and $\gamma^j$, respectively. We orient these curves by requiring that $\bar \gamma^i$ gets the orientation induced by $-\nabla(f_i-f_j)$, and $\bar \gamma^j$ gets the orientation induced by $-\nabla(f_j-f_i)$. 

Now assume that $\tilde f_i-\tilde f_j$ is another locally defined function difference associated to $\Lambda$ defined on $V\subset M$. If $U \cap V \neq \emptyset$ and if $\tilde f_l$ defines the same sheet of $\Lambda$ as $f_l$ does for $l=i,j$ on the overlap $U \cap V$, then the flow lines on $U$ being solutions to $- \nabla (f_i-f_j)$ can be extended to flow lines on $V$, now being solutions to $- \nabla (\tilde f_i-\tilde f_j)$. We denote such  patched local flow lines by \emph{flow lines of $\Lambda$}, where we are allowed to patch several times.

Thus, a flow line of $\Lambda$ is locally a flow line of the negative gradient of a locally defined function difference associated to $\Lambda$. 
The \emph{1-jet lift of a flow line of $\Lambda$} is locally given by the 1-jet lift of the corresponding local flow line. The \emph{cotangent lift of a flow line of $\Lambda$} is defined in an analogous way, and gets its orientation induced by the orientations of the cotangent lifts of the corresponding local flow lines.

\begin{figure}[ht]
\labellist
\small\hair 2pt
\pinlabel $\R$ [Br] at 95 350
\pinlabel ${f_1}$ [Br] at  400 280 
\pinlabel ${f_2}$ [Br] at 400 215
\pinlabel $f_3$ [Br] at 390 145
\pinlabel $\Gamma$ [Br] at 220 22
\pinlabel $M$ [Br] at 67 15
\endlabellist
\centering
\includegraphics[height=5cm, width=6.5cm]{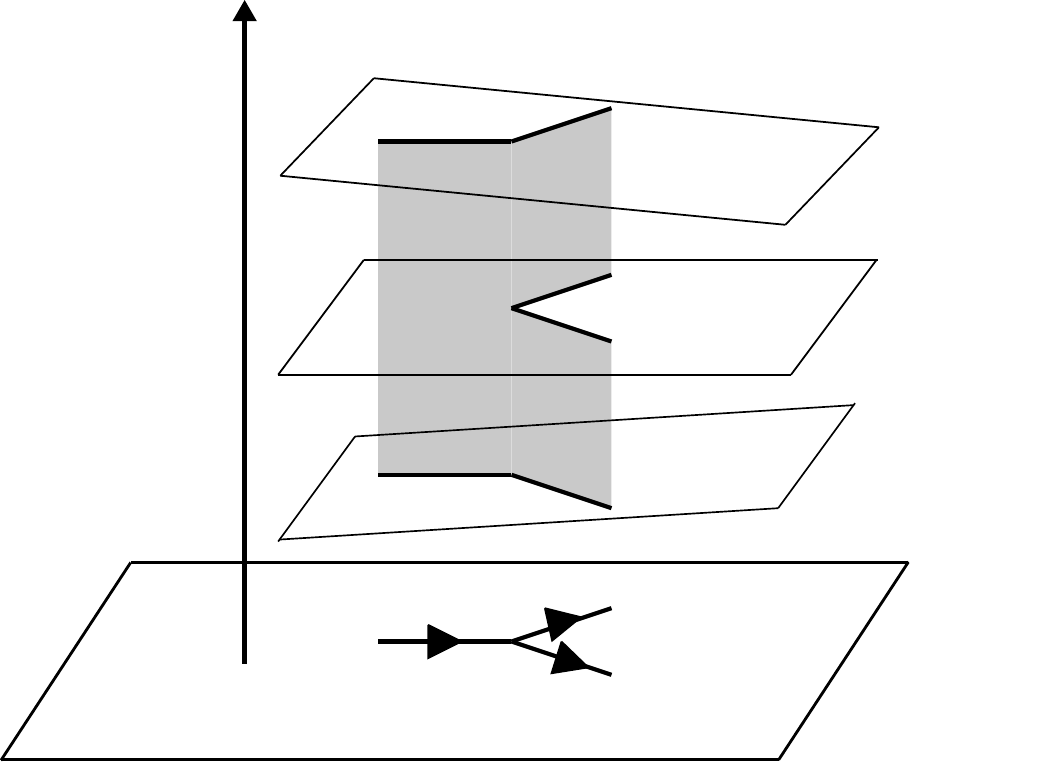}
\caption{The local picture of a Morse flow tree $\Gamma$ in a neighborhood of a $Y_0$-vertex. The graphs of the defining functions $f_1$, $f_2$ and $f_3$ are sketched, together with the lift of the tree to the sheets of $\Lambda$ determined by these functions. The shaded area indicates the corresponding pseudo-holomorphic disk.}
\label{fig:extree}
\end{figure}

The Morse flow trees are built using flow lines of $\Lambda$. That is, a \emph{Morse flow tree} $\Gamma$ is an immersion of a directed, rooted combinatorial tree into $M$ so that each edge maps to a connected subset of some flow line of $\Lambda$, and where the direction of the tree coincides with the flow direction induced by the defining vector fields of the flow line. The vertices are allowed to have valence at most 3, and at 2- and 3-valent vertices the local defining functions change. For example, at a $3$-valent vertex we may have an incoming flow line which solves $- \nabla(f_i-f_k)$, and two outgoing flow-lines solving $- \nabla (f_i-f_j)$ and $- \nabla (f_j - f_k)$, respectively. See Figure \ref{fig:extree}. 
Each edge of the tree has two lifts to $\Lambda$, given by the 1-jet lift of the  corresponding flow line. The \emph{cotangent lift of $\Gamma$} is the oriented cotangent lift of the edges of $\Gamma$ and is required to give an oriented closed curve in $\Pi_\C(\Lambda)$. The tree itself is oriented away from the root.

\begin{remark}
Let $\gamma(t)$ be a flow line of $\Lambda$. 
Note that we do not require $\gamma$ to be the solution curve of the \emph{same} function difference $- \nabla (f_i-f_j)$ everywhere, but we are allowed to switch locally defining functions along $\gamma$, without introducing 2-valent vertices. In particular, this allows $\gamma$ to be immersed, and hence also the edges of a Morse flow tree are allowed to be immersed. 
\end{remark}

The vertices of rigid Morse flow trees are of the following types:

\begin{description}
 \item[\emph{1-valent puncture}] a 1-valent vertex that is a critical point of the locally defined function difference defining the edge ending at the vertex in a neighborhood of the vertex, and thus the vertex corresponds to a Reeb chord of $\Lambda$,  
 \item[\emph{end}] 1-valent vertex whose $1$-jet lift consists of one point which is contained in $\Sigma$,
 \item[\emph{2-valent puncture}] a 2-valent vertex that lifts to Reeb chord endpoints, and thus corresponds to a Reeb chord of $\Lambda$,  
 \item[\emph{switch}] a 2-valent vertex contained in $\Pi(\Sigma)$, 
 \item[$Y_0$-vertex] a 3-valent vertex where defining functions change as in Figure \ref{fig:extree},  
 \item[$Y_1$-vertex] a 3-valent vertex similar to a $Y_0$-vertex but where the middle sheet is replaced by two sheets meeting along $\Sigma$, and the vertex belongs to $\Pi(\Sigma)$.
\end{description}

 The root of the tree is required to be a puncture, either 1- or 2-valent, and is called a \emph{positive puncture}. All other vertices that are punctures are called \emph{negative punctures}. 

For a more detailed description of the vertices  occurring in a rigid flow tree we refer to   [\cite{Georgios1}, Section 2.2] and \cite{orienttrees}.

%

\subsection{Correspondence between disks and trees}
By the results in \cite{trees}, there is a one-to-one correspondence between rigid flow trees and rigid pseudo-holomorphic disks defined by $\Lambda$. Briefly this works as follows. Consider the fiber scaling $s_\lambda:(x,y,z) \mapsto (x, \lambda y, \lambda z)$ as $\lambda \to 0$, which pushes $\Lambda$ towards the zero section. Then there is a choice of metric of $M$ and almost complex structure $J$ of $T^*M$, so that after a perturbation of $\Lambda$ there is a one-to-one correspondence between rigid flow trees of $\Lambda$ and sequences of $J$-holomorphic disks $u_\lambda$ in $T^*M$ with boundary on $\Pi_\C(s_\lambda(\Lambda))$. The disks are characterized by the property that their boundaries are arbitrarily close to the cotangent lift of the corresponding tree $\Gamma$. 
Moreover, we can identify $T^*M$ with $\C^n$ in neighborhoods of the disks, and the linearized boundary conditions of the disks tend to constant $\R^n$-boundary conditions in these coordinates, except at cusp points, where we get a split boundary condition which is constantly $\R$ in the first $n-1$ directions, and given by a uniform $\pm \pi$-rotation in the last direction, where the sign depends on the type of vertex. See Section \ref{sec:triv}. In addition, when we let $\lambda \to 0$, then  the domain of the disk will split up into a finite union of strips and strips with one slit. 

\begin{figure}[ht]
\labellist
\small\hair 1.5pt
\pinlabel $p_0$ [Br] at 45 112
\pinlabel $p_1$ [Br] at 772 8
\pinlabel $p_2$ [Br] at 772 57
\pinlabel $p_3$ [Br] at 772 106
\pinlabel $p_4$ [Br] at 772 156
\pinlabel $p_5$ [Br] at 772 203
\pinlabel $(\tau_1,1)$ [Br] at  442 26 
\pinlabel $(\tau_2,2)$ [Br] at  345 75 
\pinlabel $(\tau_3,3)$ [Br] at  485 125 
\pinlabel $(\tau_4,4)$ [Br] at  410 172 
\endlabellist
\centering
\includegraphics[height=2.5cm]{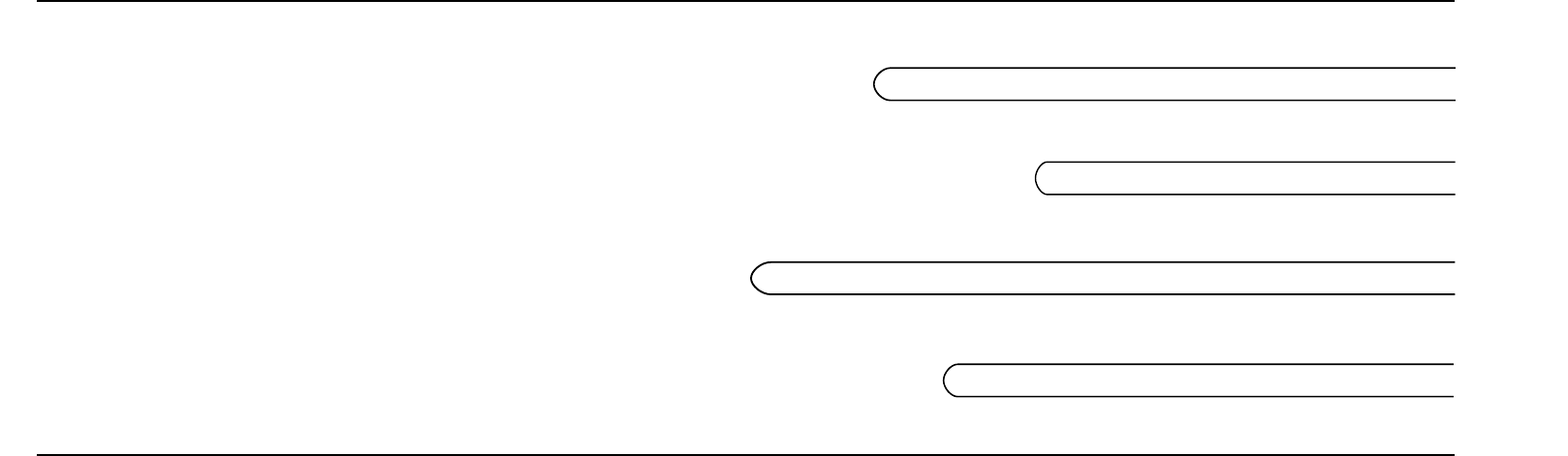}
\caption{A standard domain $\Delta_{m+1}$, $m = 5$.}
\label{fig:1}
\end{figure}

To that end, we associate a \emph{standard domain} to $\Gamma$. This is a subset $\Delta_{m+1}(\bar{\tau})\subset \R \times [0,m]\subset \R^2$,  $\bar \tau = (\tau_1,\dotsc,\tau_{m-1})$, obtained by removing $m-1$ horizontal slits starting at $(\tau_j,j)$, $j=1,\dotsc,m-1$, and ending at $\infty$. All slits look the same, they are strips ending in a half-circle of width $\epsilon$, $0< \epsilon <<1$. A point $(\tau_j,j)$ where a slit ends is called a \emph{boundary minimum}.  See Figure \ref{fig:1}. 

 The region at the negative infinity will correspond to the positive puncture of the tree, and the regions at the positive infinity to the negative punctures and the ends. The boundary minima will correspond to either 3-valent vertices or 2-valent punctures. The boundary of the standard domain should be thought of as the lifted flow lines to the cotangent lift of $\Gamma$. We write $\Delta(\Gamma)$ for the standard domain associated to $\Gamma$. The exact  value of the corresponding $\bar \tau$ will not be important in this paper, but it will be important to keep track of which boundary minimum that correspond to which vertex of the tree. We refer to [\cite{trees}, Section 2.1.A] for a more detailed discussion.

\subsection{Partial flow trees}
We will often cut a given rigid flow tree into smaller pieces, to obtain \emph{partial flow trees}. These are similar to true flow trees, except that the cotangent lift of the tree is not required to close up. The vertices which are created during the cutting are called \emph{special punctures}, and the other vertices we denote by \emph{true vertices}.

If $\Gamma$ is a rigid flow tree and $\Gamma' \subset \Gamma$ is a partial flow tree, then the edges of $\Gamma'$ inherit the directions from the edges of $\Gamma$. In particular, each partial flow tree will be rooted and oriented away from the root. If the root is a special puncture it is called a \emph{positive special puncture}. The special punctures that are not the root are called \emph{negative special punctures}. Note that when we cut a (partial) flow tree $\Gamma$ into two partial flow trees $\Gamma_1'$, $ \Gamma_2'$, then one of the trees will have the cutting point as a special positive puncture, and the other tree will have the cutting point as a special negative puncture.  Again we refer to [\cite{trees}, Section 2.2].

\begin{definition}
 A partial flow tree which is obtained by cutting a rigid flow tree just once is called a \emph{sub flow tree}. 
\end{definition}

Each partial flow tree $\Gamma'$ has an associated standard domain, which is nothing but the corresponding part of the standard domain associated to the rigid flow tree that $\Gamma'$ is a part of.

\subsection{Flow-outs and intersection manifolds}\label{sec:flowouts}
If $\gamma:I \to M$, $t \mapsto \gamma(t)$ is a flow line of $\Lambda$ that is maximally extended (that is, it cannot be extended by patching to any other  local flow line as in Section \ref{sec:trees}) we call it a \emph{maximal flow line of $\Lambda$}.
From [\cite{trees}, Lemma 2.8] we get the following.

\begin{itemize}
 \item If $I$ has a non-compact end then the flow line limit to a critical point of a locally defined function difference associated to $\Lambda$ as $t$ approaches this end. I.e.\ there is a $B \in \R$ so that $I \cap [B,\infty)= [B,\infty)$, $\lim_{t \to \infty} \gamma(t) = p$ and  $\dot \gamma|_{[B,\infty)} =-\nabla(f_1-f_2)$  for some local function difference $f_1-f_2$ which has $p$ as a critical point. (Similarly if $I \cap (-\infty,B] = (-\infty,B]$).
 \item If the maximal interval of definition has a compact end then the flow line limit to $\Pi(\Sigma)$ as $t$ approaches that end. 
\end{itemize}

This implies that if $p$ is a critical point of a locally defined function difference $f_1-f_2$ associated to $\Lambda$ we can define the stable (unstable) manifold $W^s(p)$ ($W^u(p)$) as the set of all local flow lines of  $-\nabla(f_1-f_2)$ that limit to $p$ as $t \to \infty$ ($t \to -\infty$). Note that these manifolds will only be contained in the connected set $U\subset M$ where $f_1-f_2$ is defined. 

%
%

%

There are also other interesting flow manifolds which are not necessarily related to critical points:
\begin{definition}
Let $K \subset M$ be a subset and $f_1-f_2$ a local function difference associated to $\Lambda$. Then the \emph{flow-out of K in the direction of $\pm\nabla(f_1-f_2)$} is defined as a subset $\Flo^\pm(K) \subset K \times [0,\infty)$ together with a map $i^\pm: \Flo^\pm(K) \to M$ where 
\begin{align*}
\Flo^\pm(K)=\{(p,t) \in K \times [0,\infty) & | \text{ there is a maximal flow line } \gamma_p:I \to M   \text{ of } \Lambda\\ &  
\text{ such that } I\cap[0,t]=[0,t], \gamma_p(0)=p, \\
& \text{ and there is a } B>0 \text{ such that } \dot \gamma_p|_{[0,B]} =\pm \nabla(f_1-f_2)\}
\end{align*}
and where 
\begin{equation}\label{eq:inclusion}
i^\pm: \Flo^\pm(K) \to M, \qquad (p,t) \mapsto \gamma_p(t).
\end{equation}
\end{definition}

\subsubsection{Flow outs of some particular sub flow trees}\label{sec:goodflow}
Now assume that $\Gamma' \subset \Gamma$ is a sub flow tree of a rigid tree $\Gamma$. Assume that the special puncture of $\Gamma'$ is given by $q$, and let 
 $e$ be the edge of $\Gamma'$ ending at $q$ and let $p$ denote the other vertex of $e$. Let $f_1-f_2$ be the local function differece of $\Lambda$ defining $e$ in a neighborhood of $p$. We  will define something called the \emph{flow-out of $\Gamma'$ at $q$}, denoted by $\F_q(\Gamma')$.  In the case when $p$ is a $2$- or $3$-valent vertex, then this is the flow-out of the corresponding \emph{intersection manifold $\I_p(\Gamma')$} along $e$. All this is defined inductively  below, for each type of vertex $p$.  In all cases we let $D_p^k$ be a $k$-dimensional open disk centered at $p$ and contained in the connected subset of $M$ where $f_1-f_2$ is defined. We also let $S_p \subset M$ be a sphere centered at $p$ with radius so small that it intersects $e$ in precisely one point $s$, located between $p$ and $q$.  
 
 First we assume that $e$ is a local flow line of $-\nabla{(f_1-f_2)}$ and define the flow-outs and intersection manifolds as follows.

 \begin{description}
 \item[ \it $p$  positive $1$-valent puncture] \hfill \\
  Let $\I_p(\Gamma')$ be the connected component of $S_p \cap W^u(p)$ that intersects $e$, and let $\F_q(\Gamma')$ be the flow-out of $\I_p$ in the direction of $-\nabla{(f_1-f_2)}$. 
\item [\it $p$  negative $1$-valent puncture]\hfill \\
 Let $\I_p(\Gamma')$ be the connected component of $S_p \cap W^s(p)$ that intersects $e$, and let $\F_q(\Gamma')$ be the flow-out of $\I_p$ in the direction of $\nabla{(f_1-f_2)}$. 
 \item[\it $p$ an end]\hfill \\
  Choose $\I_p(\Gamma')=D_p^{n-1} \subset \Pi(\Sigma)$ and let $\F_q(\Gamma')$ be the flow-out of $\I_p$ in the direction of $\nabla{(f_1-f_2)}$. 
  \item[\it $p$ a switch]\hfill \\
Let $r$ be a point on the edge of $\Gamma'$ containing $p$ and not containing $q$, let $p_0\neq p$ be the true vertex of $\Gamma'$ on that edge and let $\Gamma'_1$ be the sub flow tree of $\Gamma$ with $r$ as a special positive (negative) puncture if $q$ is positive (negative). Assume, by induction, that $r$ is chosen so that $\F_r(\Gamma_1')$ is defined and let $t_0>0$ so that $(p_0,t_0)\in \F_r(\Gamma'_1)$ satisfies $i(p_0,t_0) = p$, where $i:\F_r(\Gamma'_1) \to M$ is the map in \eqref{eq:inclusion}. Let $D(p_0,t_0)\subset \F_r(\Gamma'_1)$ be a neighborhood of $(p_0,t_0)$ and let  $D_p^{n-1} \subset \Pi(\Sigma)$ so that $\I_p(\Gamma')=i(D(p_0,t_0)) \cap D_p^{n-1}$ is a disk. We define $\F_q(\Gamma')$ to be the flow-out of $\I_p(\Gamma')$ in the direction of $\nabla{(f_1-f_2)}$ ($-\nabla{(f_1-f_2)}$). 
  \item[\it$p$  $2$-valent negative puncture, $q$ is positive] \hfill \\
  Define $\F_q(\Gamma')$ to be the flow-out of $p$ in the direction of $\nabla{(f_1-f_2)}$. The intersection manifold $\I_p(\Gamma')$ equals $p$. 
 \item[\it$p$ is a $Y_0$-vertex, $q$ is positive]\hfill \\
Let $e_1,e_2$ be the other edges of $\Gamma$ containing $p$, let $r_1$ and $r_2$ be points on $e_1$ and $e_2$, respectively, and let $\Gamma'_j$, $j=1,2$, be the sub flow tree of $\Gamma$ with $r_j$ as special positive puncture. Assume, by induction, that  $r_j$ is chosen so that $\F_{r_j}(\Gamma_j')$ is defined for $j=1,2$, and let $p_j \neq r_j$ be the true vertex on $e_j\cap \Gamma_j'$. Let $t_1,t_2>0$ so that $(p_j,t_j)\in \F_{r_j}(\Gamma'_j)$ satisfies $i_j(p_j,t_j) = p$, $j=1,2$, where $i_j:\F_{r_j}(\Gamma'_j) \to M$ is the map from \eqref{eq:inclusion}. Let $D(p_j,t_j)\subset \F_{r_j}(\Gamma'_j)$ be a neighborhood of $(p_j,t_j)$, $j=1,2$, and define $\I_p(\Gamma')=i_1(D(p_1,t_1)) \cap i_2(D(p_2,t_2))$.
Let $D_p^{\dim \I_p(\Gamma')} \subset  \I_p(\Gamma')$.  We define the flow-out $\F_q(\Gamma')$ to be given by the flow-out of $D_p^{\dim \I_p(\Gamma')}$  in the direction of $\nabla{(f_1-f_2)}$.      
 \item[\it $p$ is a $Y_1$-vertex, $q$ is positive]\hfill \\
 Choose $D_p^{n-1} \subset \Pi(\Sigma)$ and define $\I_p(\Gamma')=i_1(D(p_1,t_1)) \cap i_2(D(p_2,t_2))\cap D_p^{n-1}$, where $i_1(D(p_1,t_1)),i_2(D(p_2,t_2))$ are defined analogous as in the case of $Y_0$-vertices. We define the flow-out $\F_q(\Gamma')$ to be given by the flow-out of $\I_p(\Gamma')$ in the direction of $\nabla{(f_1-f_2)}$. 
\end{description}

\noindent
Next we consider the general case when the edge $e$ above is allowed to consist of patched flow lines. We assume that $q$ is not a self intersection point of $e$, and let $r$ be a point on $e$ so that $\F_r(\Gamma_1')$ is defined using the description above, where $\Gamma_1'$ is the sub flow tree of $\Gamma'$ having $r$ as special puncture. Let $t_0>0$ so that $(p,t_0)\in \F_r(\Gamma_1')$ satisfies $i(p,t_0)=q$, where $i : \F_r(\Gamma_1') \to M$ is the map in \eqref{eq:inclusion}, and let $D(p,t_0) \subset \F_r(\Gamma_1')$ be a neighborhood of $(p,t_0)$ .

\begin{definition}\label{def:tangents}
We let $\F_q(\Gamma') = \F_r(\Gamma_1')$ and we define $T_q \F_q(\Gamma')$ to be the  tangent space of $i(D(p,t_0))$ at $q$.
\end{definition}

\begin{remark}
 It might seem that the definitions of flow outs and intersection manifolds for 2--valent vertices made above are incomplete, but they are in fact sufficient for our algorithm. 
 \end{remark}

\begin{remark}\label{rmk:stablemfd}
If $p$ and $q$ belong to the same local flow line it follows that $T_q\F_q(\Gamma') \simeq T_q W^u(p)$ in the case $p$ is a positive 1-valent puncture and $T_q\F_q(\Gamma') \simeq T_q W^s(p)$ in the case $p$ is a negative 1-valent puncture .
\end{remark}

\section{The sign \texorpdfstring{$\nu_{\triv}$}{v}}\label{sec:triv}
The sign $\nu_{\triv}$ depends on the spin structure of $\Lambda$. Briefly, this has the following reason. To give a coherent orientation of the determinant line of the $\dbar$-problem corresponding to a rigid flow tree $\Gamma$, we need to choose a well-defined trivialization of $T\Lambda$ along the lift of $\Gamma$, up to homotopy. And this can be done by using the spin structure, as we will explain in this section.

\subsection{Initial choices}
To start with we need to discuss the initial orientation choices. These are the following:
 \begin{itemize}
 \item  choice of spin structure on $\Lambda$;
 \item  choice of orientation of the base manifold $M$;
 \item   choice of orientations of the  unstable manifolds associated to the critical points of the locally defining function differences for $\Lambda$;
 \item choice of orientation of  $\C$. 
\end{itemize}

 That the orientation of $\C$ affects the sign of $\Gamma$ comes from the fact that the underlying theory uses determinant line bundles associated to complex $\dbar$-operators, which are given orientations induced by the orientation of $\C$. See [\cite{fooo}, Section 8.1] and [\cite{orientbok}, Section 4.5.6 and Lemma 3.8]. This choice will affect the signs in Section \ref{sec:stab}. The choice of orientations of the unstable manifolds will be visible first in Section \ref{sec:intersection}.

\begin{remark}
 In this paper we will assume $M$ to be oriented, but in fact the moduli space of rigid flow trees of $\Lambda$ should admit a coherent orientation also in the case when $M$ in non-orientable. One can work out the algorithm for computing the signs of the trees also in this case, if one keeps track of all local choices of orientations of $M$ in the sections that follow. However, this is too technical to be covered in generality in the present article. 
\end{remark}

\subsection{Trivializations}\label{sec:trivtriv}
Let $\Gamma$ be a rigid flow tree of $\Lambda$. Using the constructions in [\cite{trees}, Section 6] we associate a model $\dbar$-problem to $\Gamma$. This model problem consists of punctured, model pseudo-holomorphic disks in $T^*M$ which are glued together so that the boundary of the model problem is arbitrarily close to the cotangent lift of $\Gamma$ as $\lambda \to 0$. The glued problem is not necessarily a solution to the $\dbar$-equation, but a Floer's Picard-like argument guarantees that one can find a true pseudo-holomorphic disk in a neighborhood of the model problem, see [\cite{trees}, Section 6.4].

The model disks that the model problem is built of capture information about the vertices and edges of $\Gamma$. Each such model problem is given an orientation and then they are glued together to give the orientation of $\Gamma$. This orientation depends on the initial orientation choices. In this section we describe how the choice of spin structure of $\Lambda$ comes into play, which has to do with a choice of trivialization of $T \Pi_\C(\Lambda)$ along the cotangent lift of $\Gamma$. This, in turn, induces a choice of trivialization of the linearized boundary condition of the linearized $\dbar$-operator corresponding to the model problem associated to $\Gamma$. We refer to [\cite{orienttrees}, Section 4] for details.

We first describe a default trivialization of $T \Pi_\C(\Lambda)$ along lifts of $\Gamma$ in neighborhoods of vertices. Then we will need to define trivializations along lifts of edges to glue these trivializations together to one that is induced by the chosen spin structure. Either one does this directly, or one uses a default trivialization also along the edges and in the end compare the default trivialization with the one induced by the spin structure.

\begin{figure}[ht]
      \labellist
\small\hair 2pt
\endlabellist
 \vspace{-11cm}
 \hspace{-2.5cm}
\centering
 \includegraphics[width=11cm, height=15cm]{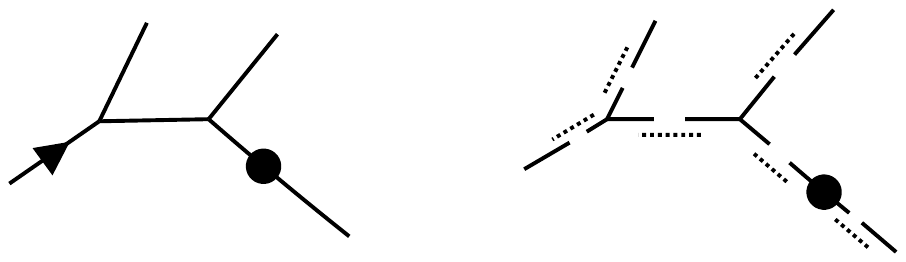} 
 \vspace{-1cm}
 \caption{A tree $\Gamma$ with one positive 1-valent puncture, 2 $Y_0$- / $Y_1$-vertices, one negative 2-valent puncture (or switch), and 3 negative 1-valent punctures /ends. To the right the corresponding subdivision into vertex regions (filled) and edge regions (dotted).}
 \label{fig:modelreg}
 \end{figure}

To that end, we divide the tree $\Gamma$ into \emph{vertex regions} and \emph{edge regions}, where vertex regions are disjoint, connected components of $\Gamma$ containing exactly one vertex each, and edge regions are complementary, connected regions, slightly overlapping the vertex regions and other edge regions but not containing any vertices. See Figure \ref{fig:modelreg}. A \emph{model region} is either a vertex region or an edge region.

The vertex regions are classified according to what type of vertex they contain, and each model region has a corresponding standard domain, as follows. 
\begin{description}
 \item[\it Positive 1-valent puncture] A standard domain $\Delta_2$ with the point at $-\infty$ representing the positive puncture.
  \item[\it Negative 1-valent puncture] A standard domain $\Delta_2$ with the point at $+\infty$ representing the negative puncture.
  \item[\it Positive 2-valent puncture] A standard domain $\Delta_3$ with the point at $-\infty$ representing the positive puncture.
  \item[\it Negative 2-valent puncture] A standard domain $\Delta_3$ with one of the points at $+\infty$ representing the negative puncture. If the negative 2-valent puncture corresponds to $p_1$ with respect to the notation in Figure \ref{fig:1}, then we say that the puncture is of \emph{Type 1}, and if it corresponds to $p_2$ we say it is of \emph{Type 2}. See Table \ref{tab:two}.
  \item[\it$Y_0$-vertex] A standard domain $\Delta_3$ with the boundary minimum representing the $Y_0$-vertex. 
    \item[\it End-vertex] A standard domain $\Delta_2$ with the point at $+\infty$ representing the end-vertex.
    \item[\it Switch] A standard domain $\Delta_2$ with a point on $\partial \Delta_2$ representing the point where the $1$-jet lift of the switch passes through $\Sigma$.
  \item[\it $Y_1$-vertex] A standard domain $\Delta_3$ with the
    boundary minimum representing the $Y_1$-vertex.
      \item[\it No true puncture] A standard domain $\Delta_2$.
\end{description}
In this way, we can interpret the standard domain $\Delta(\Gamma)$ associated to $\Gamma$ as built by the standard domains of the model regions of $\Gamma$. 

   \begin{table}[ht] \caption{{\bf Different types of negative 2-valent punctures.} }\label{tab:two}
 \begin{center}
 \tabulinesep=1.15mm
  \begin{tabu}{>{\centering\arraybackslash}m{.5in} | >{\centering\arraybackslash}m{2in}|>{\centering\arraybackslash}m{2in}}
   \hline
    Type &
  Front projection and tree &   Lagrangian projection and standard domain\\
 \hline
 1 &
     \labellist
\small\hair 2pt
\pinlabel $1$ [Br] at 218 245
\pinlabel $2$ [Br] at  218 185 
\pinlabel $3$ [Br] at 218 125
\pinlabel $p_1$ [Br] at 120 40
\pinlabel $a$ [Br] at 15 12
\pinlabel $p_2$ [Br] at 238 12
\endlabellist
\centering
\includegraphics[width=3cm, height=3cm]{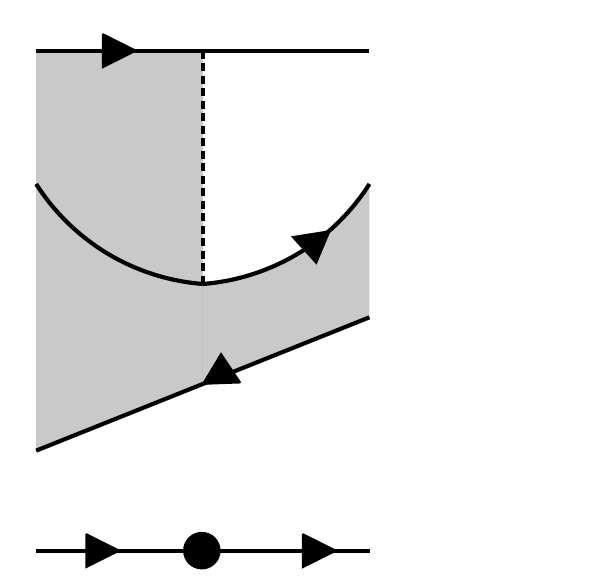} 
&
      \labellist
\small\hair 2pt
\pinlabel $1$ [Br] at 218 155
\pinlabel $2$ [Br] at  218 195 
\pinlabel $3$ [Br] at 218 235
\pinlabel $p_2$ [Br] at 228 48
\pinlabel $a$ [Br] at 30 30
\pinlabel $p_1$ [Br] at 228 5
\endlabellist
\centering
\includegraphics[width=3cm, height=3cm]{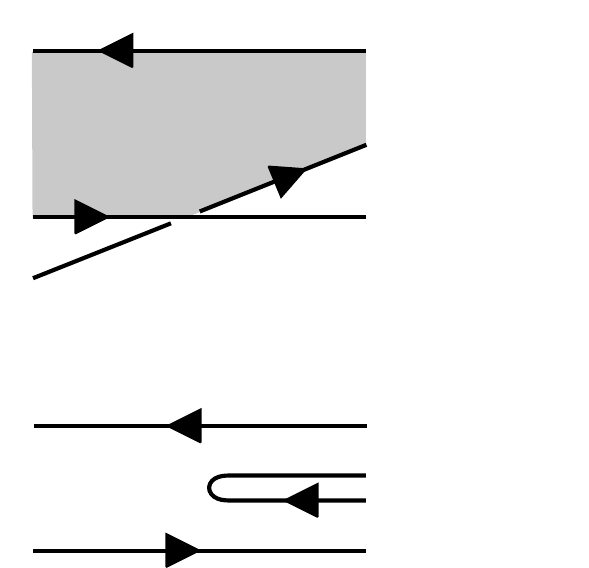} \\
\hline
  2 & 
      \labellist
\small\hair 2pt
\pinlabel $1$ [Br] at 218 190
\pinlabel $2$ [Br] at  218 120 
\pinlabel $3$ [Br] at 218 60
\pinlabel $p_1$ [Br] at 238 10
\pinlabel $a$ [Br] at 15 12
\pinlabel $p_2$ [Br] at 130 40
\endlabellist
\centering
\includegraphics[width=3cm, height=3cm]{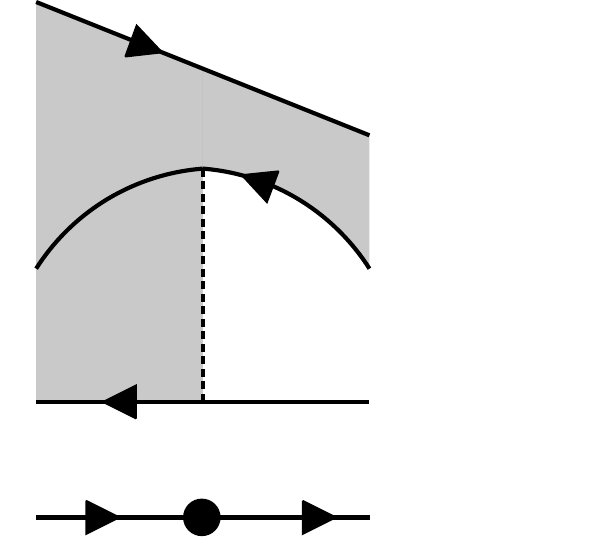} &
      \labellist
\small\hair 2pt
\pinlabel $1$ [Br] at 218 125
\pinlabel $2$ [Br] at  218 170 
\pinlabel $3$ [Br] at 218 210
\pinlabel $p_2$ [Br] at 228 53
\pinlabel $a$ [Br] at 30 30
\pinlabel $p_1$ [Br] at 228 10
\endlabellist
\centering
\includegraphics[width=3cm, height=3cm]{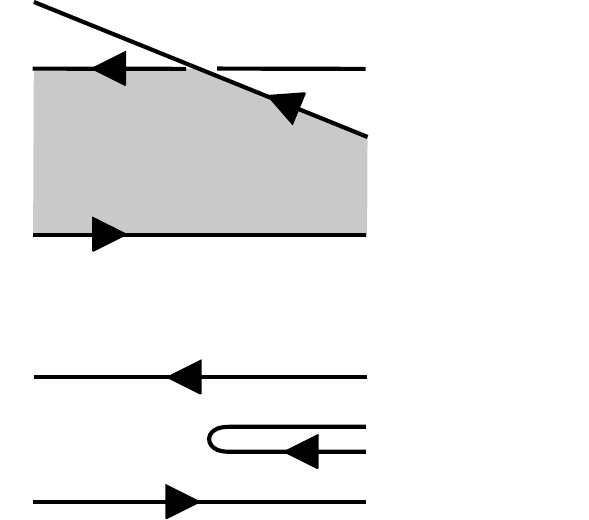} \\
   \hline
   \end{tabu}
   \end{center}

 \end{table}

Now consider the lift of such a model region $\Gamma'$, and let $\Lambda_1(\Gamma'), \dotsc,$ $\Lambda_k(\Gamma')$, $k=2,3$ or $4$ depending on the type of model region, denote the corresponding sheets of $\Lambda$. Assume that the $z$-coordinate restricted to $\Lambda_i(\Gamma')$ is greater than or equal to the $z$-coordinate restricted to $\Lambda_j(\Gamma')$ for $i < j$. Let $\tilde \gamma_i$ be the lift of $\Gamma'$ contained in $\Lambda_i(\Gamma')$, and let $\Lambda_i^\C(\Gamma')$, $\gamma_i$ be the Lagrangian projection of $\Lambda_i(\Gamma')$ and $\tilde \gamma_i$, respectively, for $1\leq i \leq k$. 

\begin{figure}[ht]
      \labellist
\small\hair 2pt
\pinlabel $\times \R^{n-1}$ [Br] at 740 142
\endlabellist
 \vspace{-4cm}
 \hspace{-2.5cm}
\centering
 \includegraphics[width=5cm, height=7cm]{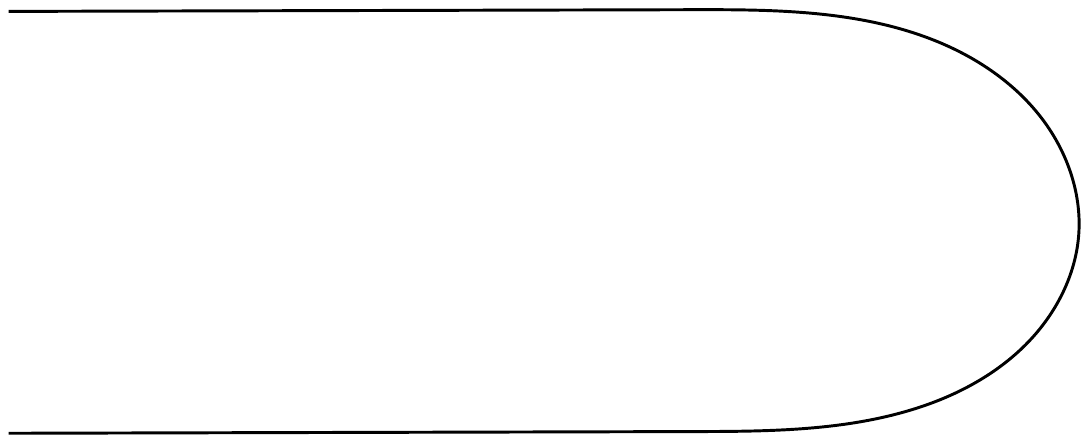} 
 \caption{Local picture in the Lagrangian projection close to an end- switch- or $Y_1$-vertex, consisting of a bended one-dimensional sheet $W$ times $\R^{n-1}$.}
 \label{fig:bend}
 \end{figure}

From the choices of metric and almost complex structure in [\cite{trees}, Section 4] together with a perturbation of $\Lambda$ to be only a totally real submanifold, also described in that section, it follows that in a neighborhood of the cotangent lift of $\Gamma'$ we can locally identify $T^*M$ with $\C^n$, where the zero section of $T^*M$ is identified with the real part $\R^n$. Moreover, using this local identification we may assume that $\Lambda_i^\C(\Gamma')$ is given by a parallel copy of $\R^n$, $1\leq i \leq k$,  in the case when $\Gamma'$ is a an edge region or a vertex region where the vertex is not contained in $\Pi(\Sigma)$. In the case when the vertex of $\Gamma'$ is contained in $\Pi(\Sigma)$ we may assume that the sheets $\Lambda_i^\C(\Gamma')$, $1 \leq i \leq k$, are given by parallel copies of $\R^n$ away from a neighborhood of $\Pi_\C(\Sigma)$, and that a pair of the sheets meet along $\Pi_\C(\Sigma)$ as in Figure \ref{fig:bend}. That is, these sheets are given by the product of $\R^{n-1}$ and a curve consisting of a half-circle joining two half-lines parallel to the remaining $\R$-factor (smoothed out).

\begin{remark}
 This gives a picture of $\Pi_\C(\Lambda)$ which is not the one we are classically used to. The reason that we might assume that the sheets in the Lagrangian projection are parallel copies of $\R^n$ away from cusps also in neighborhood of punctures comes from the fact that we are only interested in the linearized boundary condition for the $\dbar$-operator associated to the model disk, and by passing to the limit as $\lambda  \to 0$ in the same time as we put weights on the associated Sobolev spaces we can in fact consider constant $\R^n$ boundary conditions. See [\cite{orienttrees}, Sections 4 and 6].
\end{remark}

On the linearized level, this half-circle gives rise to a positive (in the end-case) or a negative (in the switch and $Y_1$-case) $\pi$-rotation in the $\C$-factor where the half-circle lives. This is defined as follows. If $V$ is a 2-dimensional real vector space spanned by orthonormal vectors $v_1, v_2$, then a positive (negative) $\pi$-rotation is given by the matrix 
\begin{equation*}
 \begin{bmatrix}
  \cos \theta & -\sin \theta \\
   \sin \theta & \cos \theta 
 \end{bmatrix},
  \qquad  
 \left(\begin{bmatrix}
  \cos \theta & \sin \theta \\
   -\sin \theta & \cos \theta 
 \end{bmatrix} \right),
 \quad \theta \in [0,\pi],
\end{equation*}
with respect to the basis $v_1, v_2$.

Using the local identification of $T^*M$ with $\C^n$ we can use the projection of  $\Lambda_i^\C(\Gamma')$ to $M$ to lift a given trivialization of $TM$ along $\Gamma'$ to a trivialization of  $T\Lambda_i^\C(\Gamma')$ along $\gamma_i$, $1 \leq i \leq k$. We will do this in a certain way, and to explain this we first need to assign a sign to the arcs $\gamma_i$, away from cusps. That is, let $\gamma \subset \gamma_i$ be a connected subarc so that  $\Lambda_i^\C(\Gamma')$ is identified with a parallel copy of $\R^n$ in a neighborhood of  $\gamma$, and let 
\begin{equation*}
 \sigma(\gamma) =
 \begin{cases}
  1, \text{ if } T\tilde \Pi: T_p\Lambda_i^\C(\Gamma') \to T_{\Pi(p)}M  \text{ is orientation-preserving for } p \in \gamma \\
    -1, \text{if } T\tilde \Pi: T_p\Lambda_i^\C(\Gamma') \to T_{\Pi(p)}M \text{ is orientation-reversing for } p \in \gamma,
 \end{cases}
\end{equation*}
where $\tilde \Pi : \Lambda_i^\C(\Gamma') \to M$ is the restriction of the base projection $\Pi: J^1(M) \to M$ to $\Lambda_i^\C(\Gamma')$.

For technical reasons we should stabilize the tangent bundle of $\Pi_\C(\Lambda)$ with an extra trivial direction, called an \emph{auxiliary direction}, and consider a trivialization of the bundle 
 \begin{equation*}                                                                                                                                                            T \tilde \Lambda := T \Pi_{\C}(\Lambda) \oplus \R,                                                                                                                                                          \end{equation*}
along the cotangent lifts of the edges of $\Gamma$.

If we also add an auxiliary direction to $TM$, to form 
\begin{equation*}
 T \tilde M = TM \oplus \R
\end{equation*}
we can canonically identify $T \tilde M$ along $\Gamma'$ with $T \tilde \Lambda$ along the cotangent lift of $\Gamma'$ away from cusps, using the identification of the sheets $ \Lambda_i^\C(\Gamma')$ with parallel copies of $\R^n$ as above. That is, if $(\partial_{x_1},\dotsc,\partial_{x_n})$ is a trivialization of $TM|_{\Gamma'} \simeq \Gamma' \times \R^n$, we extend it to a trivialization $(\partial_{x_1},\dotsc,\partial_{x_n}, \partial_{x_{n+1}})$ of $T \tilde M|_{\Gamma'}$ and lift this to a trivialization $(\partial_{x_1},\dotsc,\partial_{x_n}, \partial_{x_{n+1}})$ of $T\Lambda_i^\C(\Gamma') \oplus \R$ along $\gamma_i$, $1 \leq i \leq k$, using the inverse of the derivative of the map $p:(a_1 +x_1, \dotsc, a_n +x_n) \mapsto (x_1, \dotsc, x_n)$ for appropriate choices of $a_1, \dotsc, a_n$, and given by the identity in the auxiliary direction.

%

\begin{remark}
We assume that all trivializations of $TM$ are orientation-preserving.
\end{remark}

The trivialization of $T \tilde \Lambda$ along lifts of vertex regions $\Gamma'$ is given as follows.

 \begin{description}
 \item[\it 1-valent punctures] \hfill \\
 In this case we have two cotangent lifts $\gamma_1, \gamma_2$ of $\Gamma'$. Lift the trivialization of $T \tilde M |_{\Gamma'}$ described above to a trivialization of $T \tilde \Lambda$ along $\gamma_i$ given by
\begin{equation}\label{eq:triven}
   (\partial_{x_1},\dotsc,\partial_{x_n},\sigma(\gamma_i)\cdot\partial_{x_{n+1}}),
\end{equation}
for $i=1,2$.

  \item[\it 2-valent punctures] \hfill \\
 The trivialization is given similar to the case of 1-valent punctures. The only difference is that the cotangent lift of $\Gamma'$ now consists of three components $\gamma_1, \gamma_2, \gamma_3$. We define the trivialization along $\gamma_i$ to be given by \eqref{eq:triven} for $i =1,2,3$.

  \item[\it $Y_0$-vertices] \hfill \\
 The trivialization is given similar to the case of 2-valent punctures.
 
 \item[\it ends] \hfill \\
Let $e$ denote the end vertex. In this case we have two cotangent lifts $\gamma_1, \gamma_2$ of $\Gamma'$, which meet at $\Pi_\C(\Lambda)$. Recall that we assume $\gamma_1$ to correspond to the 1-jet lift of $\Gamma'$ with greatest $z$-coordinate. This implies that if we consider $\gamma_1, \gamma_2$ as parameterized by the corresponding parts of $\partial \Delta(\Gamma)$ then we will pass through $\Pi_\C(\Sigma)$ coming from $\gamma_1$ and leaving along $\gamma_2$, if we follow the direction induced by the orientation of $\partial \Delta(\Gamma)$. 

Assume that $T_e\Pi(\Sigma) = \spn(\partial_{x_1},\dotsc,\partial_{x_{n-1}})$, and that the union $\Lambda^\C_1(\Gamma') \cup \Lambda^\C_2(\Gamma')$ splits as $\R^{n-1} \times W \subset \R^{n-1} \times \C$, where $W$ is as in Figure \ref{fig:bend}. Let $\tilde U(e)$ be a neighborhood of the cotangent lift of $e$ so that we can find a neighborhood of $\gamma_i \setminus \tilde U(e) \subset \Lambda^\C_i(\Gamma')$ in which $\Lambda^\C_i(\Gamma')$ is a parallel copy of $\R^n$, $i =1,2$, and let $U(e) = \tilde U(e) \cap(\gamma_1 \cup \gamma_2)$. Lift the trivialization of $T \tilde M |_{\Gamma'}$ to a trivialization of $T \tilde \Lambda|_{\gamma_1 \cup \gamma_2}$ as follows. 
  \begin{itemize}
  \item Along $\gamma_i \setminus U(e)$, $i=1,2$:
  \begin{equation*}  
  (\partial_{x_1},\dotsc,\partial_{x_n},\sigma(\gamma_i\setminus U(e))\cdot\partial_{x_{n+1}})                                            \end{equation*}
\item Along  $U(e)$: Recall that we are coming in along $\gamma_1$, then we will perform a positive $\pi$-rotation in the $\C$-component corresponding to $\partial_{x_n}$ when passing the cotangent lift of $e$, where this rotation is induced by the geometry of $\Lambda$, and then we are leaving along $\gamma_2$. 
 Depending on the sign of $\sigma(\gamma_1\setminus U(e))$ we define the trivialization during the rotation as follows:
 \begin{description}
 \item[$\sigma(\gamma_1\setminus U(e))=-1$] \hfill \\
 First we perform a $-\pi$-rotation in the plane spanned by  $\partial_{x_n}$ and $\partial_{x_{n+1}}$, taking $\partial_{x_n}$ to  $-\partial_{x_n}$ and $-  \partial_{x_{n+1}}$ to $\partial_{x_{n+1}}$, keeping the other directions fixed. Assume that this rotation is made along $\gamma_1$. Then we let the trivialization be given by
 \begin{equation*}
  (\partial_{x_1},\dotsc,\partial_{x_{n-1}},-\cos(\theta) \partial_{x_n} - \sin (\theta) \partial_{y_n}, \partial_{x_{n+1}}), \quad \theta \in [0,\pi],
 \end{equation*}
 during the actual geometric rotation. Here $\partial_{y_n}$ is the vector so that the tuple $(\partial_{x_n}, \partial_{y_n})$ gives an orthonormal basis for $\C$. Assume that half of this rotation is made along $\gamma_1$, and that the other half is made along $\gamma_2$. 
 \item[$\sigma(\gamma_1\setminus U(e))=1$] \hfill \\
 In this case we begin with the geometric rotation which induces the trivialization
 \begin{equation*}
  (\partial_{x_1},\dotsc,\partial_{x_{n-1}},\cos(\theta) \partial_{x_n} + \sin (\theta) \partial_{y_n}, \partial_{x_{n+1}}), \quad \theta \in [0,\pi].
 \end{equation*}
 We assume that half of this rotation is made along $\gamma_1$, and that the other half is made along $\gamma_2$.  Just after this rotation, along $\gamma_2$, we perform a $-\pi$-rotation in the plane spanned by  $\partial_{x_n}$ and $\partial_{x_{n+1}}$, taking $- \partial_{x_n}$ to  $\partial_{x_n}$ and  $\partial_{x_{n+1}}$ to $- \partial_{x_{n+1}}$, keeping the other directions fixed.
 \end{description}
 \end{itemize}
 
 \begin{remark}
 Note that $\sigma(\gamma_1\setminus U(e)) \neq \sigma(\gamma_2\setminus U(e))$. 
 \end{remark}

 \item[\it switches] \hfill \\
 Let $s$ be the switch-vertex and note that this vertex has two cotangent lifts. Let $\tilde s$ be the lift that is contained in $\Pi_\C(\Sigma)$. We have three components $\gamma_1, \gamma_2, \gamma_3$ of the cotangent lift of $\Gamma'$, where either $\gamma_1$ and $\gamma_2$ or $\gamma_2$ and $\gamma_3$ meet at $\Pi_\C(\Sigma)$. In the case when $\gamma_1$ and $\gamma_2$ meet we pass $\tilde s$  coming from $\gamma_2$ and leaving along $\gamma_1$, and in the other case we pass  $\tilde s$ coming from $\gamma_3$ and leaving along $\gamma_2$. Consider the second case, the first case is completely similar (replace $\gamma_1$ by $\gamma_3$, $\gamma_2$ by $\gamma_1$ and $\gamma_3$ by $\gamma_2$ in the second case to get the trivialization in the first case.)
 
 Assume that $T_s\Pi(\Sigma) = \spn(\partial_{x_1},\dotsc,\partial_{x_{n-1}})$, and that the union $\Lambda^\C_2(\Gamma') \cup \Lambda^\C_3(\Gamma')$ splits as $\R^{n-1} \times W \subset \R^{n-1} \times \C$, where $W$ is as in Figure \ref{fig:bend}. Note that $\Lambda^\C_1(\Gamma')$ by assumption is identified with a parallel copy of $\R^n$.
 
 Let $\tilde U(s) \subset \Lambda^\C_2(\Gamma') \cup \Lambda^\C_3(\Gamma')$ be a neighborhood of $\tilde s$ so that we can find a neighborhood of $\gamma_i \setminus \tilde U(s) \subset \Lambda^\C_i(\Gamma')$ in which $\Lambda^\C_i(\Gamma')$ is a parallel copy of $\R^n$, $i =2,3$, and let $U(s) = \tilde U(s) \cap(\gamma_2 \cup \gamma_3)$.   Lift the trivialization of $T \tilde M |_{\Gamma'}$ to a trivialization of $T \tilde \Lambda|_{\gamma_1 \cup \gamma_2\cup \gamma_3}$ as follows. 
  \begin{itemize}
  \item Along $\gamma_i \setminus U(s)$, $i=1,2,3$:
  \begin{equation*}  
  (\partial_{x_1},\dotsc,\partial_{x_n},\sigma(\gamma_i\setminus U(s))\cdot\partial_{x_{n+1}})                                            
  \end{equation*}
\item Along  $U(s)$:
 \begin{description}
 \item[$\sigma(\gamma_3\setminus U(s))=-1$] \hfill \\
 First we perform a $-\pi$-rotation in the plane spanned by  $\partial_{x_n}$ and $\partial_{x_{n+1}}$, taking $\partial_{x_n}$ to  $-\partial_{x_n}$ and $-  \partial_{x_{n+1}}$ to $\partial_{x_{n+1}}$, keeping the other directions fixed. Assume that this rotation is made along $\gamma_3$. Then we let the trivialization be given by
 \begin{equation*}
  (\partial_{x_1},\dotsc,\partial_{x_{n-1}},-\cos(\theta) \partial_{x_n} +\sin (\theta) \partial_{y_n}, \partial_{x_{n+1}}), \quad \theta \in [0,\pi],
 \end{equation*}
 during the actual geometric rotation. Here $\partial_{y_n}$ is the vector so that the tuple $(\partial_{x_n}, \partial_{y_n})$ gives an orthonormal basis for $\C$. Assume that half of this rotation is made along $\gamma_3$, and that the other half is made along $\gamma_2$. 
 \item[$\sigma(\gamma_3\setminus U(s))=1$] \hfill \\
 In this case we begin with the geometric rotation which induces the trivialization
 \begin{equation*}
  (\partial_{x_1},\dotsc,\partial_{x_{n-1}},\cos(\theta) \partial_{x_n} - \sin (\theta) \partial_{y_n}, \partial_{x_{n+1}}), \quad \theta \in [0,\pi].
 \end{equation*}
 We assume that half of this rotation is made along $\gamma_3$, and that the other half is made along $\gamma_2$.  Just after this rotation, along $\gamma_2$, we perform a $-\pi$-rotation in the plane spanned by  $\partial_{x_n}$ and $\partial_{x_{n+1}}$, taking $- \partial_{x_n}$ to  $\partial_{x_n}$ and  $\partial_{x_{n+1}}$ to $- \partial_{x_{n+1}}$, keeping the other directions fixed.
 \end{description}
\end{itemize}

\item[\it $Y_1$-vertices] \hfill \\
This is similar to the case of switch-vertices, except that we now have two sheets $\Lambda^\C_1(\Gamma')$ and $\Lambda^\C_4(\Gamma')$ which can be identified with parallel copies of $\R^n$, and two sheets $\Lambda^\C_2(\Gamma')$ and $\Lambda^\C_3(\Gamma')$ meeting at $\Pi_\C(\Sigma)$ and whose union splits as $\R^{n-1} \times W \subset \R^{n-1} \times \C$, where $W$ is as in Figure \ref{fig:bend}. When passing through the cotangent lift of the $Y_1$-vertex that is contained $\Pi_\C(\Sigma)$, we do this by coming in along $\gamma_3$ and leaving along $
\gamma_2$. The trivialization along $\gamma_1,\gamma_2$ and $\gamma_3$ is now defined as in the case of switch-vertices, and the trivialization along $\gamma_4$ is given by \eqref{eq:triven}.

%
%
\end{description}

\begin{definition}\label{def:uplow}
If we have a switch $s$ where $\gamma_2$ and $\gamma_3$ meet at $\Pi_\C(\Sigma)$ we say that $s$ is a \emph{lower switch}, and in the other case we say that $s$ is an \emph{upper switch}. This notation will be used in Section \ref{sec:stab}.
\end{definition}

Now use the edge regions to glue together these trivializations to get a trivialization of $T\tilde \Lambda$ along the cotangent lift of the rigid flow tree $\Gamma$. A priori, one should use the spin structure of $\Lambda$ to do this, but this is not necessary as long as one compensate for the mistake one does as explained in Section \ref{sec:spinmist}. To that end, one can choose the trivialization of $T\tilde \Lambda$ along the cotanget lift of the edge regions in a way that is suitable for the geometric situation. For example, if $\Gamma'$ is an edge region with cotangent lifts $\gamma_1$ and $\gamma_2$ one could pick a trivialization $   (\partial_{x_1},\dotsc,\partial_{x_n},\partial_{x_{n+1}})$ of $T\tilde M$ along the edge regions and let the trivialization of $T\tilde \Lambda|_{\gamma_i}$ be given by $(\partial_{x_1},\dotsc,\partial_{x_n},\sigma(\gamma_i)\cdot\partial_{x_{n+1}})$, $i=1,2$. Compare with the example in Section \ref{sec:ex}. Choose a trivialization along the cotangent lift of all edge regions, use the trivializations described above for the cotangent lift of the vertex regions of $\Gamma$, and denote the resulting trivialization of $T\tilde \Lambda$ along the cotangent lift of $\Gamma$ by $\tau_d$.

%
%

\subsection{Definition of \texorpdfstring{$\nu_{\triv}$}{v}}\label{sec:spinmist}
Our chosen trivialization $\tau_d$ should now be compared with a trivialization induced by the spin structure of $\Lambda$. 

In fact, the choice of spin structure gives a trivialization of $T \tilde \Lambda$ over the $1$-skeleton of $\Lambda$ that extends to the $2$-skeleton. Assuming that the $0$-skeleton is chosen so that all Reeb chord end points are contained in it, homotope the trivialization induced by the spin structure over the 1-skeleton to agree with $\tau_d$ at all Reeb chord end points.  Then homotope the lifts of the edges of $\Gamma$, keeping the homotopy constant along the $0$-skeleton, to the $1$-skeleton and let the spin structure induce a trivialization, which we denote by $\tau_s$. The fact that the spin structure extends to the $2$-skeleton implies that the trivialization we get is independent of the chosen homotopy (up to homotopy, which is everything we need as far as orientations are concerned). See  [\cite{orientbok}, Section 3.4.2].

How this is done explicitly depends on the specific case, and is difficult to describe in full generality. We explain how this is done for one-dimensional Legendrian knots in Section \ref{sec:ex}.

To make the dependence of spin structure, and hence the dependence of the Lagrangian trivialization along the cotangent lift of $\Gamma$, to something that can be put into our algorithm, we proceed as follows. For each connected component $l$ of the cotangent lift of $\Gamma$ we compare the trivialization $\tau_d$ and $\tau_s$ along $l$. By our assumptions the transition maps from $\tau_d$ to $\tau_s$ along $l$ induce a loop in $SO(n+1)$. Using that  $\pi_1(SO(n+1)) \simeq \Z_2$, $n \geq 1$, we define $\nu(l) = 0$ if this loop is homotopically trivial, and $\nu(l) = 1$ otherwise.

 \begin{definition}\label{def:triv}
  Let 
\begin{equation*}
L = \{\text{connected components of the cotangent lift of } \Gamma\}
\end{equation*}
and let 
\begin{equation}\label{eq:beta}
 \nu_{\triv}(\Gamma) = \prod_{l \in L} (-1)^{\nu(l)}.
\end{equation}
\end{definition}

%

\section{The sign \texorpdfstring{$\nu_{\inter}$}{v}}\label{sec:intersection}
In this section we will define the sign $\nu_{\inter}(\Gamma)$. This is nothing but a geometric intersection sign, coming from intersections of the flow-outs of the sub flow trees of $\Gamma$. Indeed, from [\cite{trees}, Proposition 3.14] it follows that the tangent space, in the space of trees, of a sub flow tree $\Gamma'$ with special puncture $q$ can be identified with $T_q\F_q(\Gamma')$,  using the notation from Definition \ref{def:tangents}, via evaluation at the puncture. In particular, the dimension of the space of sub flow trees is never bigger than the dimension of $M$. 

We use this observation to define orientations of the sub flow trees of $\Gamma$ given in Section \ref{sec:goodflow}, realized as orientations of $T_q\F_q(\Gamma')$. This will be done inductively, over the number of true vertices of the sub flow trees, by taking oriented intersections of the flow-outs in $M$. In the end, we will be able to define an orientation of the tangent space of $\Gamma$, and since this is $0$-dimensional it is only a sign. This will be the sign $\nu_{\inter}(\Gamma)$.   


\subsection{Intersection orientation}
To calculate the orientations of the intersection manifolds of $\Gamma$ we will make use of the following orientation rule.

Assume that $V_1, V_2 \subset \R^n$ are  subspaces so that the sequence 
  \begin{equation}
 \begin{array}{ccccccccc}\label{eq:glu102}
0 
&\to 
& V
&\to
& 
V_1  \oplus
V_2
&\to 
& \R^n
&\to 
& 0 \\
& &v &\mapsto &(v,v) &   & & &\\
& & & & (u,w) &\mapsto & w-u,  & & 
\end{array}
\end{equation} 
is exact, where $V = V_1 \cap V_2$. Assume that $V_1$, $V_2$ and $\R^n$ are oriented, denote the chosen orientation by $\Or(V_1), \Or(V_2), \Or(\R^n)$, respectively. If $\dim V_1 > \dim V$ or $\dim V_2 > \dim V$, pick an  orientation $\Or_c(V)$ of $V$ and let $\tilde V_i \subset V_i$, $i=1,2$, be oriented subspaces so that 
\begin{equation}\label{eq:vent}
 V_i = V \oplus \tilde V_i \quad \text{ oriented}, \quad i=1,2.
\end{equation}
Then let $\nu \in \{0,1\}$ so that 
\begin{equation}\label{eq:intor}
\bigwedge^{\max} V \wedge \bigwedge^{\max} \tilde V_1 \wedge \bigwedge^{\max} \tilde V_2 = (-1)^\nu \bigwedge^{\max} \R^n \quad \text{ oriented.}
\end{equation}

If $\dim V = \dim V_1 = \dim V_2$, let  $\Or_c(V) =  \Or(\R^n)$ and let  $\nu \in \{0,1\}$ so that 
\begin{equation*}
 \Or(V_1) = (-1)^\nu\Or(V_2).
\end{equation*}

\begin{definition}
We define the \emph{intersection orientation} $\Or_{\inter}(V)$ of $V$ to be given by
\begin{equation*}
 \Or_{\inter}(V) := \Or_{\inter}(V_1, V_2)=(-1)^{\nu+\dim V_1\cdot(1+\dim V)}\Or_c(V).
\end{equation*}

\end{definition}

In what follows, we will use the orientation of $M$ to identify $T_qM \simeq \R^n$, oriented, for any $q \in \Gamma$.

\subsection{Orientation of flow-outs and intersection manifolds}\label{sec:flowor}
Now we define the orientations of the sub flow trees of $\Gamma$, using the inductively construction of the flow-outs from Section \ref{sec:goodflow}. We use the notation from there, so that $q$ is a special puncture of a sub flow tree $\Gamma'$ of $\Gamma$ and $e$ is the edge of $\Gamma'$ between $q$ and the vertex $p$. Here we do not require $e$ to be a local flow line (that is, we allow it to be patched). Note that, in all cases except for when $p$ is a switch we have that $\F_q(\Gamma') \simeq \I_p \times [0,\infty)$ (or $\F_q(\Gamma') \simeq \I_p \times [0,A]$ for some $A > 0$). In the case when $p$ is a switch we have that $\F_q(\Gamma') \simeq \I_p$. Hence, from Remark \ref{rmk:stablemfd} it follows that if we fix an orientation of $W^u(p)$ ($W^s(p)$) in the case when $p$ is a positive (negative) puncture we get an induced orientation $\Or (T_q(\F_q(\Gamma'))=:\Or(\F_q(\Gamma'))$ of $\F_q(\Gamma')$. Similarly, if we fix  an orientation of $T_pM$ in the case when $p$ is an end, $\I_p$ in the case when $p$ is a switch, and $\I_p \times e$ in the case when $p$ is 3-valent or a 2-valent puncture, we also get an induced orientation $\Or(T_q(\F_q(\Gamma'))=\Or(\F_q(\Gamma'))$ in these cases. Now we will make a choice of all these orientations.

\begin{figure}[ht]
\labellist
\small\hair 2pt
\pinlabel $\Pi^{-1}(v)$ [Br] at  35 234 
\pinlabel ${x_1}$ [Br] at 55 80
\pinlabel $x_2$ [Br] at 3 45
\pinlabel $\Pi(\Lambda_v)$ [Br] at 125 15
\pinlabel $M$ [Br] at -50 15
\pinlabel $\Lambda$ [Br] at 100 160
\endlabellist
\centering
\includegraphics[height=5cm]{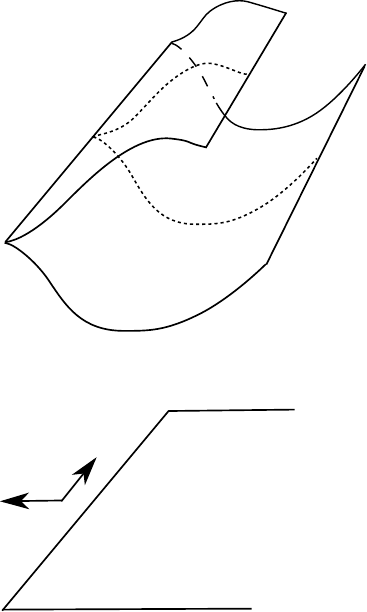}
\caption{Local coordinates at a cusp. Here $v$ represents a switch, end or $Y_1$-vertex, and $\Pi(\Lambda_v)$ is given the orientation from $\partial_{x_1}$.}
\label{fig:cusp}
\end{figure}

First we let $\Or_{\capp}(W^u(p))$ denote the initial choice of orientation of $W^u(p)$, represented by a wedge product of an oriented basis of $T_p W^u(p)$, and let it induce an orientation $\Or_{\capp}(W^s(p))$ on $W^s(p)$ by requiring that
\begin{align*}
\bigwedge^{\max} T_pM = \Or_{\capp}(W^u(p)) \wedge \Or_{\capp}(W^s(p))
\end{align*} 
is an oriented identification. Also, recall that $\Sigma \subset \Lambda$ denotes the singular set under the base projection $\Pi: \Lambda \to M$ and that $v \in \Pi(\Sigma)$ whenever $v$ is an end-, switch- or $Y_1$-vertex. At such points we may assume that $T_v\Pi(\Sigma)$ is of codimension 1 in $T_vM$ and we can locally consider $\Pi(\Sigma)$ as the boundary of $\Pi(\Lambda_v)$, where $\Lambda_v$ is (one of) the sheets of $\Lambda$ containing the lift of $v$ that is contained in $\Sigma$. Give $\Pi(\Lambda_v)$ the orientation from $M$, and define an orientation $\Or_{\capp}(T_v\Pi(\Sigma))$ to be given by the orientation of $\Pi(\Sigma)$ oriented as the boundary of $\Pi(\Lambda_v)$, oriented by outward normal last. See Figure \ref{fig:cusp}.

Now let the flow-out orientation of the  first four trees in Section \ref{sec:goodflow} be given as follows:

\noindent
{\tabulinesep=1.2mm
\begin{tabu}{{>{\arraybackslash}m{1.9in}  >{\arraybackslash}m{3.9in}}}
{\it $p$ pos.\ $1$-valent puncture:}
& 
 $\Or(\F_q(\Gamma'))= \Or_{\capp}(W^u(p))$\\
{\it $p$ neg.\ $1$-valent puncture:}&
$\Or(\F_q(\Gamma'))=\Or_{\capp}( W^s(p))$\\
{\it $p$ an end:} &
  $\Or(\F_q(\Gamma'))= \Or(T_qM)$, $n \geq 2$ \\
  & $\Or(\F_q(\Gamma'))= \Or(e)$, $n =1$, where $e$ is the edge of $\Gamma$ ending at $p$ oriented against the flow direction. 
\\
 {\it $p$  a switch, $q$ positive:}&  
 $ \Or(\F_q(\Gamma'))=\Or_{\inter}(\Or(T_p\F_p(\Gamma'_1)) ,\Or_{\capp}(T_p\Pi(\Sigma))) $ \\
  {\it $p$  a switch, $q$ negative:}&  
 $ \Or(\F_q(\Gamma'))=\Or_{\inter}(\Or_{\capp}(T_p\Pi(\Sigma)),\Or(T_p\F_p(\Gamma'_1))) $
   \end{tabu}}
   
   \noindent
   where $\Gamma'_1$ is the sub flow tree of $\Gamma$ having $p$ as positive (negative) special puncture in the case $q$ is positive (negative).

For the other trees from Section \ref{sec:goodflow}, with positive special punctures, we let the orientation of the intersection manifolds be defined as follows.
\begin{description}
 \item[ \it$p$ a negative $2$-valent puncture]   
\begin{equation*}
   \Or(T_p\I_p(\Gamma'))=\Or_{\inter}(\Or(T_p\F_p(\Gamma'_1)), \Or_{\capp}(T_pW^s(p))),
 \end{equation*} 
   where $\Gamma'_1$ is the sub flow tree of $\Gamma$ having $p$ as positive special puncture, 
 \item[\it$p$ a $Y_0$-vertex, $q$ is positive] 
 \begin{equation*}
   \Or(T_p\I_p(\Gamma'))=\Or_{\inter}(\Or(T_p\F_p(\Gamma_1')),\Or(T_p\F_p(\Gamma_2'))),
  \end{equation*}   
 \item[\it$p$ a $Y_1$-vertex, $q$ is positive]  
 \begin{equation*}
 \Or(T_p\I_p(\Gamma'))=\Or_{\inter}\left(\Or(T_p\F_p(\Gamma_1')),\Or_{\inter}(\Or(T_p\F_p(\Gamma_2')), \Or_{\capp}(T_p\Pi(\Sigma)))\right)
 \end{equation*} 
\end{description}
where in the two latter cases $\Gamma_1'$ and $\Gamma_2'$ are the sub flow trees of $\Gamma$ having $p$ as special positive puncture and where they are numbered so that the standard domain of $\Gamma_1'$ corresponds to the lower part of the standard domain of $\Gamma'$.

To define the flow-out orientation of the three last cases, let $f_1 > f_2$ be the defining functions for $e$ in a neighborhood of $p$, let $w = \nabla (f_1 - f_2)$ at $p$ (that is, $w$ is the vector tangent to $e$ at $p$ pointing in the direction against the defining flow), and define the orientation of the flow-out along $e$ as
  \begin{equation*}
  \Or(\F_q(\Gamma'))=\Or(T_p\I_p(\Gamma')) \wedge w.
  \end{equation*}

\subsection{Definition of \texorpdfstring{$\nu_{\inter}$}{v}}\label{sec:finglu}
To finally define $\nu_{\inter}(\Gamma)$ we need to handle two separate cases, depending on the valence of the positive puncture $a$ of $\Gamma$.

\begin{figure}[ht]
\labellist
\small\hair 2pt
\pinlabel ${p_2}$ [Br] at  198 400 
\pinlabel ${p_1}$ [Br] at  406 400 
\pinlabel ${a}$ [Br] at  285 470 
\pinlabel ${a}$ [Br] at 1005 470
\pinlabel $\Gamma_{2}'$ [Br] at 115 480
\pinlabel $\Gamma_a$ [Br] at 305 390
\pinlabel $\Gamma_1'$ [Br] at 485 480
\pinlabel $\Gamma$ [Br] at 1135 480
\pinlabel $q_6$ [Br] at 40 515
\pinlabel $q_7$ [Br] at 40 370
\pinlabel $q_5$ [Br] at 610 570
\pinlabel $q_4$ [Br] at 690 555
\pinlabel $q_3$ [Br] at 699 460
\pinlabel $q_2$ [Br] at 645 390
\pinlabel $q_1$ [Br] at 615 335
\pinlabel $a$ [Br] at 830 130
\pinlabel $q_1$ [Br] at 1310 7
\pinlabel $q_2$ [Br] at 1310 53
\pinlabel $q_3$ [Br] at 1310 98
\pinlabel $q_4$ [Br] at 1310 142
\pinlabel $q_5$ [Br] at 1310 185
\pinlabel $q_6$ [Br] at 1310 228
\pinlabel $q_7$ [Br] at 1310 273
\endlabellist
\centering
\includegraphics[height=5.cm]{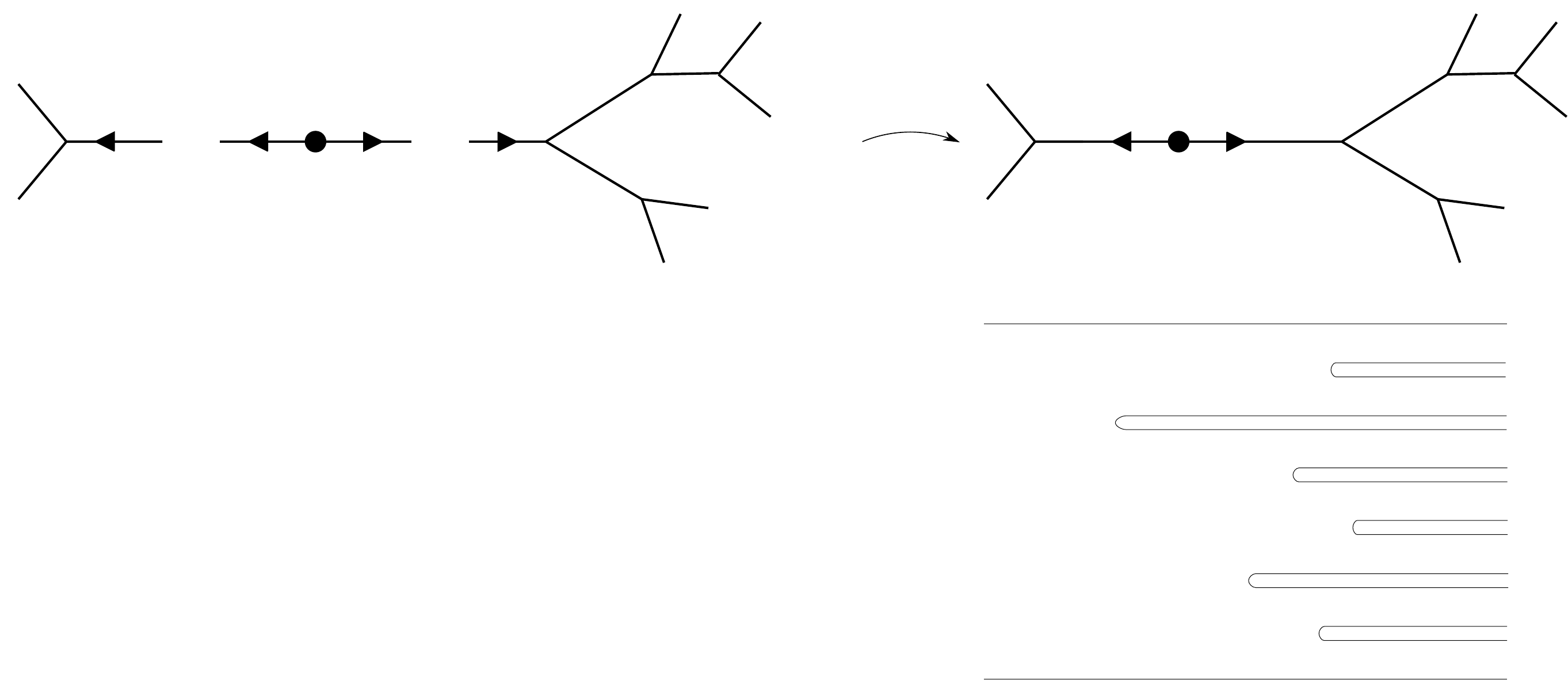}
\caption{The final gluing when the positive puncture is $2$-valent. 
}
\label{fig:finish2val}
\end{figure}

\begin{description}
\item[Case 1, $a$ is 1-valent] 
We cut $\Gamma$ into two sub flow trees, as follows.

\begin{itemize}
 \item Let $v$ denote the first vertex of $\Gamma$ that we meet when going along $\Gamma$, starting at the positive puncture $a$, satisfying that $v$ is not a  switch. 
 \item Let $e$ be the edge of $\Gamma$ ending at $v$, with respect of the orientation of the tree.
\item Let $p$ denote a point on $e$, let $\Gamma_1'$ ($\Gamma_2'$) be the sub flow tree of $\Gamma$ having $p$ as special positive (negative) puncture.
 \end{itemize}
 
Now we get two different cases, depending on the number of true vertices of $\Gamma_1'$. 
 \begin{description}
  \item[Case 1a), $\Gamma'_1$ has only one true vertex] 
If $\Gamma_1'$ has only one true vertex $v$, then assume that  the flow direction of the edge $e$  is given by $w$ at $p$.
Let $\mu_{\inter}$ be an integer so that the intersection orientation of $T_{p}\F_{p}(\Gamma_1') \cap T_{p} \F_{p}(\Gamma_2')$ is given by 
 \begin{equation}\label{eq:nu11}
  \Or_{\inter}(\Or(T_p\F_{p}(\Gamma_1')), \Or(T_p\F_{p}(\Gamma_2'))) = (-1)^{\mu_{\inter}}w.
 \end{equation}
 
 \item[Case 1b), $\Gamma_1'$ has more than one true puncture] 
 Recall that we have \linebreak inductively defined orientations of $T_v\F_{v}(\Gamma_2')$ and $T_v\I_v(\Gamma_1')$. Let $\mu_{\inter}$ be an integer so that 
\begin{equation}\label{eq:nu12}
 \bigwedge^{\max}T_vM =(-1)^{\mu_{\inter}} \cdot \Or(T_v\I_v(\Gamma_1')) \wedge \Or(T_v\F_{v}(\Gamma_2')). 
\end{equation}

\begin{figure}[ht]
\labellist
\small\hair 2pt
\pinlabel ${a}$ [Br] at  35 385 
\pinlabel ${p}$ [Br] at 310 380
\pinlabel $s$ [Br] at 170  360
\pinlabel $\Gamma_{2}'$ [Br] at 100 280
\pinlabel $\Sigma$ [Br] at 275 310
\pinlabel $\Gamma_1'$ [Br] at 450 280
\pinlabel $\Gamma$ [Br] at 1140 280
\pinlabel Glue [Br] at 760 390
\pinlabel $v$ [Br] at 445 380
\pinlabel $a$ [Br] at 850 100
\pinlabel $p_3$ [Br] at 580 488
\pinlabel $p_1$ [Br] at 650 480
\pinlabel $p_2$ [Br] at 670 370
\pinlabel $p_5$ [Br] at 620 295
\pinlabel $p_4$ [Br] at 580 245
\pinlabel $p_1$ [Br] at 1343 10
\pinlabel $p_2$ [Br] at 1343 57
\pinlabel $p_3$ [Br] at 1343 100
\pinlabel $p_4$ [Br] at 1343 143
\pinlabel $p_5$ [Br] at 1343 189
\endlabellist
\centering
\includegraphics[height=4.5cm]{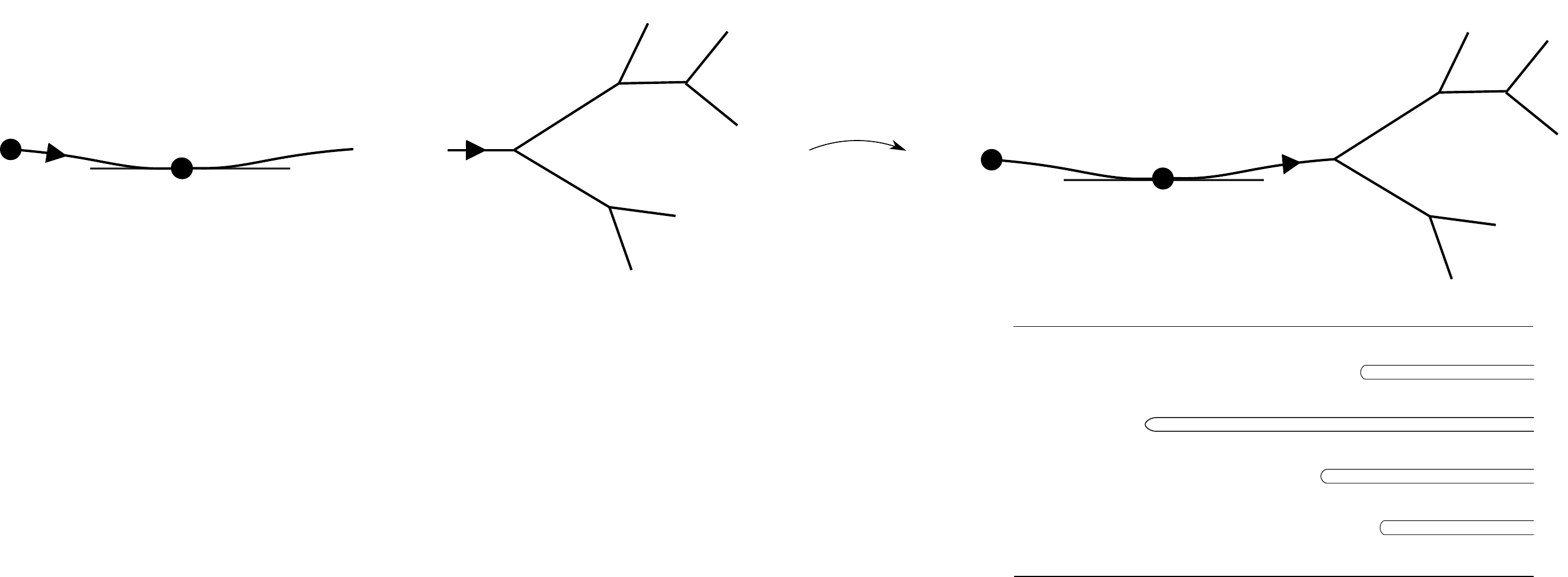}
\caption{The final gluing when the positive puncture $a$ is 1-valent, giving the rigid flow tree. }
\label{fig:finish1val}
\end{figure}

\begin{remark}
Note that we are using the intersection manifold of $\Gamma_1'$ in the second case, instead of the flow-out. 
\end{remark}
\end{description}
\item[Case 2), $a$ is $2$-valent] Let $p_1$, $p_2$ be points on each of the edges of $\Gamma$ that contain $a$. By cutting $\Gamma$ at these points we get the following three partial flow trees.
\begin{itemize}
 \item A partial flow tree $\Gamma_a$ that contains the positive 2-valent puncture $a$, and has $p_1$ and $p_2$ as negative special punctures, and no other vertices.
 \item Sub flow trees $\Gamma_i'$ with $p_i$ as a positive special puncture, $i=1,2$.  
\end{itemize}
We use the following notation, see Figure \ref{fig:finish2val}. 
\begin{itemize}
 \item Let $\Gamma_1'$ denote the sub flow tree that corresponds to the lower part of the standard domain of $\Gamma$.
 \item Let $\mu_0 \in\{0,1\}$ so that 
 \begin{equation*}
 \Or_{\capp}(W^u(a)) = (-1)^{\mu_0}.
 \end{equation*}
\item Let $\mu_i \in \{0,1\}$ so that 
\begin{equation*}
  \Or(T_{p_i}M) = (-1)^{\mu_i}\Or(T_{p_i}\F_{p_i}(\Gamma_i)), \qquad i=1,2.
 \end{equation*}
\end{itemize}
Finally let
\begin{equation}\label{eq:nu13}
\mu_{\inter} = \mu_0 + \mu_1 + \mu_2.
\end{equation}

\begin{remark}
Notice that in this case we cannot have a switch adjacent to the positive puncture. 
\end{remark}

\end{description}

\vspace{.5cm}

\begin{definition}\label{def:nui}
The sign $\nu_{\inter}$ is given by
\begin{equation*}
\nu_{\inter}(\Gamma) = (-1)^{\mu_{\inter}},
\end{equation*}
where $\mu_{\inter}$ is given by 
\begin{itemize}
\item \eqref{eq:nu11} in the case when the standard domain of $\Gamma$ has no slits,
\item  \eqref{eq:nu12} when the standard domain of $\Gamma$ has at least one slit and the positive puncture of $\Gamma$ is $1$-valent,
\item \eqref{eq:nu13} in the case when the positive puncture of $\Gamma$ is 2-valent.
\end{itemize}

\end{definition}

 \section{The sign \texorpdfstring{$\nu_{\en}$}{v}}\label{sec:end}
 Now we discuss the sign $\nu_{\en}$, coming from regarding ends as marked points of the standard domains. 
 
 Assume that $\Gamma$ has $l$ negative punctures $q_1,\dotsc,q_l$ and $k$ end-vertices $e_1,\dotsc,e_
 k$. Pick indices $i_1<\dotso<i_k, j_1<\dotso<j_l$ so that 
 \begin{equation*}
  \{i_1,\dotsc,i_k,j_1,\dotsc,j_l\} = \{1,\dotsc,k+l\}
 \end{equation*}
and so that the ends and punctures of $\Gamma$ are occurring along the boundary of $\Delta(\Gamma)$ in a way so that  $p_{j_r}=q_r$, $r=1,\dotsc,l$, and $p_{i_r} =e_r$, $r=1,\dotsc,k$, with notation as in Figure \ref{fig:1}. Let $m =k+l$.
 
For $m\geq 2$ we identify the ordered tuple $(p_1,\dotsc,p_m)$ with the oriented standard basis for $\R^{m}$, so that we can use wedge products of the points when computing $\nu_{\en}$. 

\begin{figure}[ht]
\labellist
\small\hair 1.5pt
\pinlabel $q_0$ [Br] at 25 112
\pinlabel $e_1$ [Br] at 762 8
\pinlabel $q_1$ [Br] at 762 57
\pinlabel $q_2$ [Br] at 762 106
\pinlabel $e_2$ [Br] at 762 156
\pinlabel $q_3$ [Br] at 762 203
\endlabellist
\centering
\includegraphics[height=2.5cm]{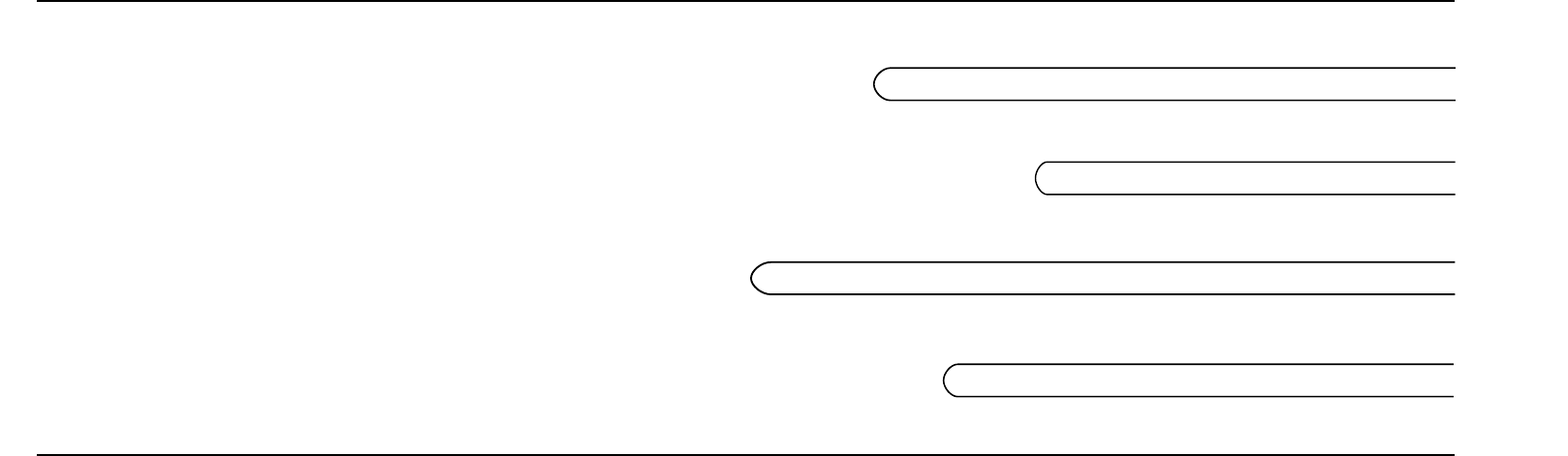}
\caption{
In this example we have $j_1=2$, $j_2=3$, $j_3=5$, $i_1=1$, $i_2=4$ and $\nu_{\en}=1$, since $e_1 \wedge q_{1} \wedge q_{2} \wedge e_2 \wedge q_3 = q_1 \wedge q_2 \wedge q_3 \wedge e_1 \wedge e_2$.}
\label{fig:endconf}
\end{figure}

 \begin{definition}\label{def:musign}
 The sign $\nu_{\en}$ is given by $(-1)^{\mu_{\en}}$, where $\mu_{\en}$ satisfies
%
%
  \begin{align*}
  (-1)^{\mu_{\en}}  p_1  \wedge  p_{2} \wedge  \dotsm \wedge  p_m 
    =  p_{j_1} \wedge p_{j_2} \wedge \dotsm \wedge p_{j_l} \wedge p_{i_1} \wedge p_{i_2} \wedge \dotsm \wedge p_{i_k}.
   \end{align*}
 See Figure \ref{fig:endconf}.
 \end{definition}

\section{The sign \texorpdfstring{$\nu_{\stab}$}{v}}\label{sec:stab}
In this section we define the algebraic sign $\nu_{\stab}$, which comes from something called capping orientations of the pseudo-holomorphic disks corresponding to the trees, and from gluing of such disks.
 Again we refer to \cite{orienttrees} for the underlying theory.

This sign will be stated as a function depending on the vertices of $\Gamma$. In particular, each vertex of $\Gamma$ will have a sign associated to it, and the product of all these signs gives $\nu_{\stab}$. 
To that end, let $\Pu_1^-$ denote the set of negative $1$-valent punctures of $\Gamma$, $\Pu_2^-$ the set of negative $2$-valent punctures, $\Y_i$ the set of $Y_i$-vertices, $i=0,1$, $\Sw$ denote the set of switch-vertices and $\Pu^+$ the set of positive punctures. The set $\Sw$ we divide into two subsets $\Sw^-$ and $\Sw^+$, respectively, where $\Sw^+$ is the set of switches that belong to the same component as the positive puncture in the partition of $\Gamma$ in Section \ref{sec:finglu}. We also let $\Pu^- = \Pu^-_1 \cup \Pu^-_2$.

Now, for each of these six sets we define a function
\begin{equation}\label{eq:vertfunc}
\sigma_{\U}: \U \to \Z, \qquad \U= \Pu_1^-,\Pu_2^-, \Y_0, \Y_1, \Sw, \Pu^+.
\end{equation}
These functions will depend on the combinatorics of $\Gamma$ in an inductive way, and they also depend on the geometry of $\Lambda$ in a neighborhood of the vertex. 

The dependence on the geometry of $\Lambda$ can be reduced to a local dependence as proven in \cite{orienttrees}, and in that paper it is proven that this only depends on the \emph{intersection parity}, the \emph{Morse index} and the \emph{sheet signs} of the vertices of $\Gamma$, where this data is defined as follows.

First, if $p$ is a Reeb chord of $\Gamma$ we define the \emph{intersection parity} $|\mu(p)|$ to be the modulo 2 count of intersections between $\Sigma$ and a path in $\Lambda$ connecting the start and end point of $p$, assuming the path intersects $\Sigma$ transversely.

\begin{definition}
If $p$ is a puncture of a Morse flow tree $\Gamma$ we let $|\mu(p)|$ be the \emph{intersection parity} of the corresponding Reeb chord. If $p$ is a special puncture which is not contained in $\Pi(\Sigma)$ we define $|\mu(p)|$ in a completely analogous way.
\end{definition}

\begin{definition}
If $p$ is a Reeb chord of $\Lambda$, then we let $I(p)$ denote the Morse index of $p$, i.e.
\begin{equation*}
 I(p) = \dim W^u(p).
 \end{equation*} 
 \end{definition}

 Note that, close to a Reeb chord, the tangent map of the restriction $\tilde \Pi$ of the base projection  $\Pi: J^1(M) \to M$ to $\Lambda$ induces isomorphisms  $T\tilde \Pi: T_{p_\pm}\Lambda \to T_{\Pi(p)}M$, where $p_{\pm}$ is the 1-jet lift of $p$ and where we choose notation so that $z(p_+)>z(p_-)$. 

 \begin{definition}
We define the \emph{sheet signs} of $p$ to be the tuple $(k_+,k_-) \in \{0,1\} \times \{0,1\}$ such that  
\begin{equation}\label{eq:sheetsign}
 \begin{cases}
  k_i=0, \text{ if } T \tilde \Pi: T_{p_i}\Lambda \to T_{\Pi(p)}M  \text{ is orientation-preserving}, \\
  k_i=1, \text{if } T\tilde \Pi: T_{p_i}\Lambda \to T_{\Pi(p)}M \text{ is orientation-reversing},
 \end{cases}
\end{equation}
$i \in \{+,-\}$. 
\end{definition}

If $p$ is a special puncture which is not contained in $\Pi(\Sigma)$ we define the sheet signs of $p$ in an analogous way.

\begin{remark}
 If $|\mu(p)| = 0$, then the sheet signs of $p$ is either $(0,0)$ or $(1,1)$, and if $|\mu(p)| = 1$, then the sheet signs of $p$ is either $(0,1)$ or $(1,0)$.
\end{remark}

\begin{definition}
 If $\Gamma$ is a (partial) flow tree, then we let $e(\Gamma)$ be the number of end vertices of $\Gamma$.
\end{definition}

We will also need to use data from the sub flow tree(s) that end(s) at $p$ as in the description in Section \ref{sec:goodflow}. In particular, if $\Gamma'$ is a (sub) flow tree we let $\bm(\Gamma')$ denote the number of boundary minima of the corresponding standard domain. For trees $\Gamma'$ with $\bm(\Gamma')>0$ we moreover let $\bm_{\min}(\Gamma')$ be the order of the boundary minimum of $\Gamma'$ with smallest $\tau$-value.  For example, if $\Gamma'$ is obtained from gluing two trees $\Gamma_1'$ and $\Gamma_2'$ at a $Y_0$- or $Y_1$-vertex, where the trees are ordered as explained in Section \ref{sec:flowor}, then 
\begin{equation*}
\bm_{\min}(\Gamma')= \bm(\Gamma_1') +1, \qquad \bm(\Gamma')-\bm_{\min}(\Gamma')= \bm(\Gamma_2').
\end{equation*}
If $\bm(\Gamma')=0$ we let $\bm_{\min}(\Gamma')=0$.

\begin{example}
 In Figure \ref{fig:finish1val} we have $\bm_{\min}(\Gamma)= 3$, and in Figure \ref{fig:finish2val} we  have $\bm_{\min}(\Gamma)= 5$. 
\end{example}

We will also use a certain subspace of the flow-out of $\Gamma'$ at the special puncture $q$, called the \emph{true kernel} of $\Gamma'$, $\Ker_q(\Gamma')$. 
This space is related to the kernel of the corresponding linearized $\dbar$-operator.
We will omit the $q$ in the notation of $\Ker_q(\Gamma')$  when there is no risk of confusion.

The space $\Ker(\Gamma')$ is defined inductively as follows, with the set-up introduced in Section \ref{sec:goodflow}. In particular, $q$ is the special puncture of $\Gamma'$ and $p$ denotes a true vertex of $\Gamma'$ connected to the special puncture $q$ via the edge $e$. However, here we do not assume $e$ to be a local flow line of $\Lambda$; we are allowed to change defining function difference along $e$, but $q$ is not allowed to be a point where $e$ intersect itself.

\begin{description}
\item[\it $p$ pos.\ $1$-valent puncture] \hfill \\
Using the notation in Definition \ref{def:tangents} we define $\Ker_q(\Gamma') = i(D(p,t_0))$ (which implies that $T_q \Ker_q(\Gamma') = T_q \F_q(\Gamma')$).
\item[\it $p$ neq.\ $1$-valent puncture] \hfill \\
Using the notation in Definition \ref{def:tangents} we define $\Ker_q(\Gamma') = i(D(p,t_0))$ (which implies that $T_q \Ker_q(\Gamma') = T_q \F_q(\Gamma')$).
\item[\it $p$ end] \hfill \\
Using the notation in Definition \ref{def:tangents} we define $\Ker_q(\Gamma') = i(D(p,t_0))$ (which implies that $T_q \Ker_q(\Gamma') = T_q M$).
\item[ \it $p$ neg.\ $2$-valent puncture]\hfill \\
$\Ker_q(\Gamma') = \emptyset$.
\item[\it $p$  a switch] \hfill \\
Assume that $\Gamma'$ has $q$ as special positive (negative) puncture. Let $\Gamma'_1 \subset \Gamma'$ be the sub flow tree with $p$ as special positive (negative) puncture and assume by induction that $\Ker_p(\Gamma'_1)$ is defined. 
If $\Ker_p(\Gamma'_1) = \emptyset$ we define $\Ker_q(\Gamma') = \emptyset$. Otherwise, let $D^{n-1}_p \subset \Pi(\Sigma)$ (where we use the notation from Section \ref{sec:goodflow} where $D^k_p$ means a $k$-dimensional disk centered at $p$), let $\pk_q(\Gamma')$ be the flow-out of $\Ker_p(\Gamma'_1) \cap D^{n-1}_p$ in the direction of $(-)\nabla(f_1-f_2)$, where $f_1-f_2$ is the defining function difference for $e$ in a neighborhood of $p$. Let $(p,t_0) \subset \pk_q(\Gamma')$ so that $i((p,t_0)) = q$, where $i:\pk_q(\Gamma') \to M$ is the map in \eqref{eq:inclusion} and let $D(p,t_0) \subset \pk_q(\Gamma')$ be a neighborhood of $(p,t_0)$ so that $i(D(p,t_0))$ is a disk. Now we get two different cases:
\begin{description}
\item[$T_pe \in T_p\Ker_p(\Gamma'_1)$]\hfill \\
In this case we define $\Ker_q(\Gamma') = i(D(p,t_0))$.
\item[$T_pe \notin T_p\Ker_p(\Gamma'_1)$]\hfill \\
Let $D_q^{n-1}$ be transverse to $e$ at $q$ so that $i(D(p,t_0))\cap D_q^{n-1}$ is a disk and let $\Ker_q(\Gamma') = i(D(p,t_0))\cap D_q^{n-1}$.
\end{description}
\item[\it $p$ a $Y_0$-vertex, $q$ is positive]\hfill \\
Let $\Gamma_1', \Gamma_2' \subset \Gamma'$ be the sub flow trees with $p$ as special positive puncture and assume by induction that $\Ker_p (\Gamma'_i)$ is defined for $i=1,2$. If $\Ker_p (\Gamma'_1) =\emptyset$ or $\Ker_p (\Gamma'_2)=\emptyset$ we define $\Ker_q(\Gamma') = \emptyset$. Otherwise, let $\pk_q(\Gamma')$ be the flow-out of $\Ker_p(\Gamma'_1) \cap \Ker_p(\Gamma'_2)$ in the direction of $\nabla (f_1-f_2)$, where $f_1-f_2$ is the locally defining function difference for $e$ in a neighborhood of $p$. Let $(p,t_0) \subset \pk_q(\Gamma')$ so that $i((p,t_0)) = q$, where $i:\pk_q(\Gamma') \to M$ is the map in \eqref{eq:inclusion} and let $D(p,t_0) \subset \pk_q(\Gamma')$ be a neighborhood of $(p,t_0)$. Let $D_q^{n-1}$ be transverse to $e$ at $q$ and so that $i(D(p,t_0))\cap D_q^{n-1}$ is a disk, and let $\Ker_q(\Gamma') = i(D(p,t_0))\cap D_q^{n-1}$.
\item[\it $p$ a $Y_1$-vertex, $q$ is positive]\hfill \\
Let $\Gamma_1', \Gamma_2' \subset \Gamma'$ be the sub flow trees with $p$ as special positive puncture and assume by induction that $\Ker_p (\Gamma'_i)$ is defined for $i=1,2$. If $\Ker_p (\Gamma'_1) =\emptyset$ or $\Ker_p (\Gamma'_2)=\emptyset$ we define $\Ker_q(\Gamma') = \emptyset$. Otherwise, choose $D_p^{n-1} \subset \Pi(\Sigma)$ and let 
 $\pk_q(\Gamma')$ be the flow-out of $D_p^{n-1}\cap \Ker_p(\Gamma'_1) \cap \Ker_p(\Gamma'_2)$ in the direction of $\nabla (f_1-f_2)$, where $f_1-f_2$ is the locally defining function difference for $e$ in a neighborhood of $p$. Let $(p,t_0) \subset \pk_q(\Gamma')$ so that $i((p,t_0))  = q$, where $i:\pk_q(\Gamma') \to M$ is the map in \eqref{eq:inclusion} and let $D(p,t_0) \subset \pk_q(\Gamma')$ be a neighborhood of $(p,t_0)$. Let $D_q^{n-1}$ be transverse to $e$ at $q$ so that $i(D(p,t_0))\cap D_q^{n-1}$ is a disk and let $\Ker_q(\Gamma') = i(D(p,t_0))\cap D_q^{n-1}$.
\end{description}

%

We are now ready to describe the functions in \eqref{eq:vertfunc}.

\subsection{Definition of \texorpdfstring{$\sigma_{\Pu_1^-}$}{s}}

Let $\Reeb$ denote the set of all Reeb chords of $\Lambda$. Then $\sigma_{\Pu_1^-}$ is given as a function
\begin{align*}
\sigma_{\Pu_1^-}: \Reeb &\to \{0,1\}\\ 
p & \mapsto  \sigma_{\Pu_1^-}(p)=\sigma_{\Pu_1^-}((k,k'),I(p)).
\end{align*}
where $(k,k')$ is the sheet signs of $p$.

If $n>1$ and $I(p)>0$ we have
\begin{align*}
 \sigma_{\Pu_1^-}((l,l),I(p)) &= n\cdot (I(p) + 1), \qquad l \in \{0,1\},\\
 \sigma_{\Pu_1^-}((0,1),I(p)) &= n\cdot I(p) +I(p) + 1, \\
 \sigma_{\Pu_1^-}((1,0),I(p)) &=  n\cdot I(p) +I(p) + 1, 
\end{align*} 
and if $n>1$, $I(p)=0$ we have
\begin{align*}
 \sigma_{\Pu_1^-}((l,l),0) &=  1,\qquad l \in \{0,1\},\\
   \sigma_{\Pu_1^-}((0,1),0) &= n, \\
\sigma_{\Pu_1^-}((1,0),0) &= n.
\end{align*} 

For $n=1$ we get the following relations:
\begin{align*}
 \sigma_{\Pu_1^-}((0,0),1) &=  0, \\
  \sigma_{\Pu_1^-}((1,1),1) &= 1,\\
 \sigma_{\Pu_1^-}((0,1),1) &= \sigma_{1,-}^{(0,1)}, \\
 \sigma_{\Pu_1^-}((1,0),1) &=  1  \\
 \sigma_{\Pu_1^-}((l,l),0) &=  1 + l + \sigma_{0,+}^{(l,l)} ,\qquad l \in \{0,1\},\\
   \sigma_{\Pu_1^-}((1,0),0) &= \sigma_{0,+}^{(0,0)} + 1, \\
\sigma_{\Pu_1^-}((0,1),0) &= \sigma_{0,+}^{(1,1)}+\sigma_{1,-}^{(0,1)}.
\end{align*} 

Here $ \sigma_{1,-}^{(0,1)},\sigma_{0,+}^{(0,0)}, \sigma_{0,+}^{(1,1)}$ are some auxiliary signs that are introduced in [\cite{orienttrees}, Section 9] and which depend on the chosen orientation of $\C$. 

\subsection{Definition of \texorpdfstring{$\sigma_{\Y_0}$}{Y}}
\begin{figure}[ht]
      \labellist
\small\hair 2pt
\pinlabel $p$ [Br] at 60 100
\pinlabel $p_1$ [Br] at  135 190 
\pinlabel $p_2$ [Br] at 130 70
\pinlabel $\Gamma_1'$ [Br] at 250 151
\pinlabel $\Gamma_2'$ [Br] at 240 70
\endlabellist
\centering
\includegraphics[height=2.5cm]{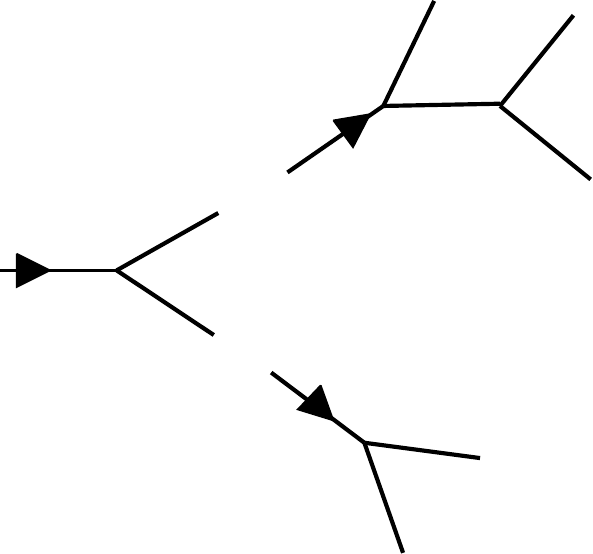} 
\hspace{1.8cm}
\caption{Gluing of two trees at a $Y_0$- or $Y_1$-vertex.}
\label{fig:gluedy0flow}
\end{figure}

Let $p$ be a $Y_0$-vertex and let $\Gamma_1'$, $\Gamma_2'$ be the sub flow trees of $\Gamma$ with positive special punctures at $p$. Here we are using the notation from Section \ref{sec:flowor}, and in particular the same ordering of the trees as in there. Consider the model region of $p$ and let $p_0$ be the corresponding positive special puncture and $p_1$, $p_2$ the negative punctures.
By cutting the trees $\Gamma_1'$ and $\Gamma_2'$ a bit later than $p$ (in the flow direction), we can as well assume that $\Gamma_i'$ has positive special puncture given by $p_i$, $i=1,2$. See Figure \ref{fig:gluedy0flow}.

 We now give the formula for the sign $\sigma_{\Y_0}(p)$, which will depend on data from the incoming trees $\Gamma_1'$, $\Gamma_2'$, the sheet signs of $p_0$, and also on the intersection parity of the special punctures $p_1$ and $p_2$.
 
That is, there is a function $\sigma_0:\{0,1\}^4 \to \{0,1\}$ so that 
we have
\begin{align*}
\sigma_{\Y_0}(p)&=\sigma_0((k,k'),|\mu(p_1)|,|\mu(p_2)|) + \eta \\
&\quad{} + \bm_{\min}(\Gamma_1') + \bm_{\min}(\Gamma_2')+\bm(\Gamma_1') + \bm(\Gamma_2'),
\end{align*}
where 
\begin{align*}
\eta&=  e(\Gamma_1')\cdot[\bm(\Gamma'_2) + e(\Gamma'_2) +1] +n \cdot \dim \F_{p_2}(\Gamma_2') + \dim \Ker_{p_2} (\Gamma_2') \\
&\quad{} + \bm(\Gamma_1')\cdot[|\mu(p_2)| + \dim \F_{p_2}(\Gamma_2') + n +1] + \dim \Ker_{p_0}(\Gamma)  +\dim \Ker_{p_1} (\Gamma_1')\\
&\quad{}+ [n+\dim \F_{p_1}(\Gamma_1')+|\mu(p_1)|]\cdot[1 + |\mu(p_2)| + \bm(\Gamma'_2)],
\end{align*}
\begin{align*}
&\sigma_0((1, 1),0,0)  \equiv \sigma_0((0, 0),0,0)   \equiv n+1, \\
& \sigma_0((0, 1),0,1)  \equiv \sigma_0((1, 0),0,1)  \equiv n+1, \\ 
& \sigma_0((0, 1),1,0)  \equiv \sigma_0((1, 0),1,0)    \equiv n+1, \\
& \sigma_0((0, 0),1,1)  \equiv \sigma_0((1, 1),1,1)  \equiv n, \pmod{2}, 
\end{align*}
and where $(k,k')$ is the sheet signs of $p_0$.

\subsection{Definition of \texorpdfstring{$\sigma_{\Y_1}$}{Y}}
With completely the same notation as in the case of $Y_0$-vertices we define the sign $\sigma_{\Y_1}$ as follows.

There is a function $\sigma_1:\{0,1\}^4 \to \{0,1\} $ so that 
we have
\begin{align*}
\sigma_{\Y_1}(p)&= 
\sigma_1((k,k'),|\mu(p_1)|,|\mu(p_2)|) + \eta +  \dim \F_{p_1}(\Gamma_1') +  \dim \F_{p_2}(\Gamma_2')+1\\ &\quad{}+\bm_{\min}(\Gamma_1') +  \bm_{\min}(\Gamma_2')+\bm(\Gamma_1') + \bm(\Gamma_2'),
\end{align*}
where, if $n>1$,
\begin{align*}
&\sigma_1((0,0),1,0) \equiv  \sigma_1((1,1),1,0) \equiv  \sigma_1((1,0),1,1) \equiv  \sigma_1((0,1),1,1)   \equiv n + 1,  \\ 
&\sigma_1((0,0),0,1) \equiv  \sigma_1((1,1),0,1)  \equiv n, \\
 &\sigma_1((0,1),0,0)  \equiv  \sigma_1((1,0),0,0)  \equiv  n+1, \pmod{2},
\end{align*}
 and where $(k,k')$ is the sheet signs of $p_0$.
 
 If $n=1$ we have the following relations:
\begin{align*}
&\sigma_1((0,0),1,0) \equiv  \sigma_1((1,1),1,0) \equiv  \sigma_1((1,0),1,1) \equiv  \sigma_1((0,1),1,1)   \equiv 1,  \\ 
&\sigma_1((0,0),0,1) \equiv  \sigma_1((1,1),0,1)  \equiv \sigma_1((0,1),0,0)  \equiv  \sigma_1((1,0),0,0)  \equiv  0, \pmod{2}.
 \end{align*}

\subsection{Definition of \texorpdfstring{$\sigma_{\Pu^-_2}$}{s}}
\begin{figure}[ht]
      \labellist
\small\hair 2pt
\pinlabel $p_1$ [Br] at  218 92 
\pinlabel $p$ [Br] at 110 65
\pinlabel $\Gamma_1'$ [Br] at 300 130
\endlabellist
\centering
\includegraphics[height=2.5cm]{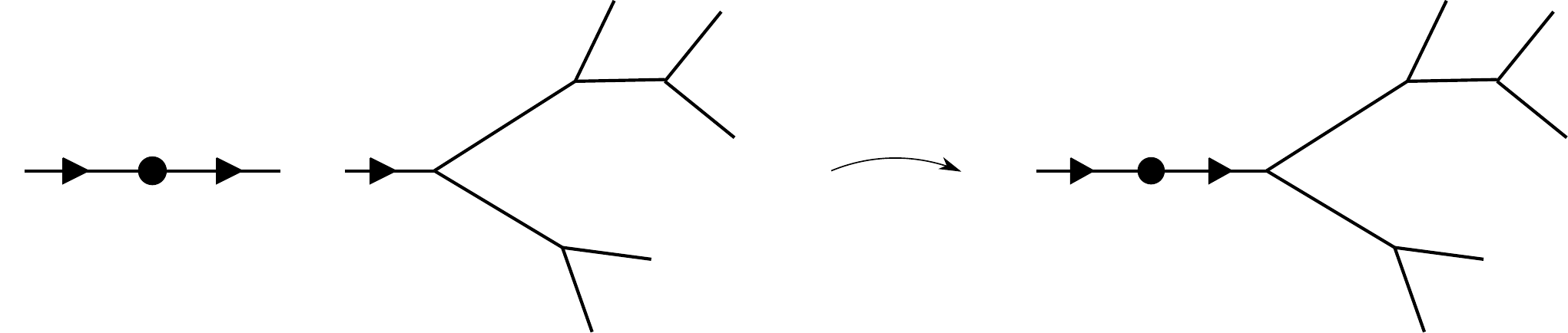} 
\hspace{1.8cm}
\caption{Gluing a tree $\Gamma'_1$ to a negative 2-valent vertex.}
\label{fig:glue2whole}
\end{figure}

In the case when $p$ is a negative 2-valent vertex we have one sub flow tree $\Gamma'_1$ which has $p$ as a special positive puncture, and again we can cut $\Gamma'_1$ a bit later, at a special puncture $p_1$, say. See Figure \ref{fig:glue2whole}.
Recall from Section \ref{sec:trivtriv}, Table \ref{tab:two}, that there are two different types of negative 2-valent punctures.
We let $\tp(p) \in\{1,2\}$ indicate the type of $p$ as given in that table. 

Just as in the case for $1$-valent negative punctures we have a function 
\begin{align*}
\sigma_{\negg,2}: \Pu_2^- &\to \{0,1\}, \\
p &\mapsto \sigma_{\negg,2}(p)= \sigma_{\negg,2}((k,k'),|\mu(p)|, |\mu(p_1)|, \tp(p))
\end{align*}
where $(k,k')$ is the sheet sign for the positive puncture on the model region containing $p$.

The sign $\sigma_{\Pu^-_2}$ is defined as follows.

\begin{description}
 \item[$p$ of Type 1]
 \begin{align*}
\sigma_{\Pu^-_2}(p) &= \sigma_{\negg,2}(p)+
\bm_{\min}(\Gamma_1')+\bm(\Gamma_1') + \dim \Ker_{p_1}(\Gamma_1')\\
&\quad{} +[n +|\mu(p)|]\cdot [\bm(\Gamma_1') +  1+ |\mu(p_1)|],
\end{align*}
 \item[$p$ of Type 2]
\begin{align*}
\sigma_{\Pu^-_2}(p) &= \sigma_{\negg,2}(p)+
\bm_{\min}(\Gamma_1')+\bm(\Gamma_1') + \dim \Ker_{p_1}(\Gamma_1')\\
&\quad{} + e(\Gamma_1') +\bm(\Gamma_1')\cdot[n + |\mu(p)| + 1].
\end{align*}
\end{description}

When $n>1$ we moreover have that 
\begin{align*}
 \sigma_{\negg,2}((0,0),0, 0, 1) &\equiv  \sigma_{\negg,2}((0,0),0, 0, 2) \equiv   \sigma_{\negg,2}((0,1),0, 1, 1) \equiv n+1\\
 &\equiv   \sigma_{\negg,2}((1,0),0, 1, 2) + 1, \\
  \sigma_{\negg,2}((1,1),0, 0, 1) &\equiv   \sigma_{\negg,2}((1,1),0, 0, 2) \equiv   \sigma_{\negg,2}((1,0),0, 1, 1)  \equiv n+1\\
  &\equiv   \sigma_{\negg,2}((0,1),0, 1, 2)+ 1, \\
  \sigma_{\negg,2}((0,1),1, 0, 1) &\equiv  \sigma_{\negg,2}((0,1),1, 0, 2) \equiv
  \sigma_{\negg,2}((0,0),1, 1, 1) \equiv   \sigma_{\negg,2}((1,1),1, 1, 2) \equiv n, \\
    \sigma_{\negg,2}((1,0),1, 0, 1) &\equiv  \sigma_{\negg,2}((1,0),1, 0, 2) \equiv
 \sigma_{\negg,2}((1,1),1, 1, 1) \equiv  \sigma_{\negg,2}((0,0),1, 1, 2) \equiv n.\\
  \end{align*}
  
   For $n=1$ we have the following relations:
  \begin{align*}
 \sigma_{\negg,2}((0,0),0, 0, 1) &\equiv  \sigma_{\negg,2}((0,1),0, 1, 1) \equiv \sigma_{0,+}^{(0,0)}, \\
  \sigma_{\negg,2}((1,0),0, 1, 2) &\equiv   \sigma_{\negg,2}((0,0),0, 0, 2)+ 1 \equiv \sigma_{0,+}^{(0,0)}+1, \\
  \sigma_{\negg,2}((1,1),0, 0, 2) &\equiv  \sigma_{\negg,2}((0,1),0, 1, 2) +1 \equiv \sigma_{0,+}^{(1,1)} + 1, \\
    \sigma_{\negg,2}((1,0),0, 1, 1) &\equiv  \sigma_{\negg,2}((1,1),0, 0, 1) \equiv \sigma_{0,+}^{(1,1)} + 1, 
  \end{align*}
  
    \begin{align*}
      \sigma_{\negg,2}((1,0),1,0,1) &\equiv \sigma_{0,+}^{(0,0)} + 1, \\
      \sigma_{\negg,2}((1,0),1,0,2) &\equiv \sigma_{0,+}^{(1,1)} + 1 ,\\
     \sigma_{\negg,2}((0,1),1,0,1) &\equiv  \sigma_{0,+}^{(1,1)} + \sigma^{(0,1)}_{1,-} ,\\
     \sigma_{\negg,2}((0,1),1,0,2) &\equiv \sigma_{0,+}^{(0,0)} + \sigma^{(0,1)}_{1,-},\\
 \end{align*}

  \begin{align*}
\sigma_{\negg,2}((0,0),1, 1, 1) +  \sigma_{\negg,2}((1,1),1, 1, 2) & \equiv \sigma_{0,+}^{(1,1)} +  \sigma_{0,+}^{(0,0)},\\
\sigma_{\negg,2}((0,0),1, 1, 1) +  \sigma_{\negg,2}((0,0),1, 1, 2) & \equiv
\sigma_{0,-}^{(0,1)} +  1,  \\
 \sigma_{\negg,2}((1,1),1, 1, 1) +  \sigma_{\negg,2}((1,1),1, 1, 2)& \equiv \sigma_{0,-}^{(0,1)} +  1,  \\
\sigma_{\negg,2}((1,1),1, 1, 1) +  \sigma_{\negg,2}((0,0),1, 1, 2) & \equiv \sigma_{0,+}^{(1,1)} +  \sigma_{0,+}^{(0,0)}.
  \end{align*}

\subsection{Definition of \texorpdfstring{$\sigma_{\Sw^-}$}{s}}\label{sec:switchsign}
\begin{figure}[ht]
      \labellist
\small\hair 2pt
\pinlabel $p_1$ [Br] at  365 97 
\pinlabel $p$ [Br] at 170 50
\pinlabel $\Gamma_1'$ [Br] at 450 150
\endlabellist
\centering
 \includegraphics[height=1.9cm]{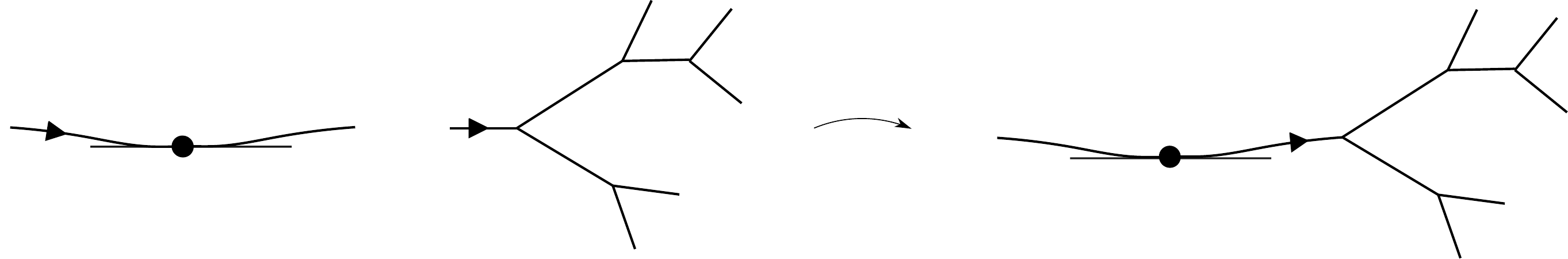} 
 \hspace{1.8cm}
 \caption{Gluing a tree $\Gamma'_1$ to a switch.}
 \label{fig:glueswichp}
 \end{figure}

Similar to the case of negative 2-valent vertices we have one sub flow tree $\Gamma'_1$ which has $p$ as a special positive puncture, and again we can cut $\Gamma'_1$ a bit later, at a special puncture $p_1$. See Figure \ref{fig:glueswichp}. The sign of $p$ is defined as follows.

First, we have a function

\begin{align*}
\sigma_{\sw}: \Sw &\to \{0,1\}, \\
p &\mapsto \sigma_{\sw}(p)= \sigma_{\sw}((k,k'),\tp(p))
\end{align*}
where $(k,k')$ is the sheet signs of the positive puncture $p_0$ of the model region of $p$, and $\tp(p)=l$ if  $p$ is a lower switch and $\tp(p)=u$ if $p$ is an upper switch, as defined in Definition \ref{def:uplow}. This function is given as follows
\begin{align*}
 &\sigma_{\sw}((0, 0),l) \equiv \sigma_{\sw}((1, 1),l) \equiv 0, \\
 &\sigma_{\sw}((0, 0),u) \equiv  \sigma_{\sw}((1, 1),u) \equiv 1, \\
 &\sigma_{\sw}((1, 0),u) \equiv  \sigma_{\sw}((0, 1),u) \equiv \sigma_{\sw}((1, 0),l) \equiv  \sigma_{\sw}((0, 1),l)  \equiv 0,
\end{align*}
everything modulo 2.

The sign of the switch $p$ is then defined as 
\begin{equation*}
\sigma_{\Sw^-}(p) = \sigma_{\sw}(p) + 1+\dim \Ker_{p_1}(\Gamma_1') +\dim \Ker_{p_0}(\Gamma) +  \dim \F_{p_1}(\Gamma_1').
\end{equation*}

\subsection{The sign of the positive puncture}
We use the notation from Section \ref{sec:finglu}. In particular, the sign of the positive puncture $a$ will depend on whether it is $1$- or $2$-valent. 

Similar to the case of 2-valent negative punctures there is a function
\begin{align*}
\sigma_{\pos,2}: \Pu_2^+ &\to \{0,1\}, \\
p &\mapsto \sigma_{\pos,2}(p)= \sigma_{\pos,2}((k,k'),|\mu(q_1)|, |\mu(q_2)|)
\end{align*}
where $p_1, p_2$ are the special punctures as described in Section \ref{sec:finglu}, $(k,k')$ is the sheet signs of $p$ and where $\Pu_2^+$ is the set of all possible positive $2$-valent punctures of trees of $\Lambda$. This function is given by \begin{equation*}
 \sigma_{\pos,2}((0,0),0,0) = \sigma_{\pos,2}((1,1),0,0)=\sigma_{\pos,2}((0,1),1,0)=\sigma_{\pos,2}((1,0),1,0)=0                                                                                                                                                                                                                                                                                                                       \end{equation*}
 for all $n$.

   For $n>1$: 
  \begin{align*}
   &\sigma_{\pos,2}((0,1),0,1) \equiv \sigma_{\pos,2}((1,0),0,1) \equiv 0, \\
   &\sigma_{\pos,2}((0,0),1,1)  \equiv \sigma_{\pos,2}((1,1),1,1)  \equiv 1, \pmod{2}. 
 \end{align*}
 
 For $n=1$:
  \begin{align*}
   &\sigma_{\pos,2}((0,1),0,1) \equiv \sigma_{\pos,2}((1,0),0,1) \equiv \sigma_{0,+}^{(0,0)} + \sigma_{0,+}^{(1,1)}. \\
  \end{align*}

Let $k$ denote the number of elements in $\Sw^+$, and note that $k=0$ if $a$ is $2$-valent. Then the sign $\sigma_{\Pu^+}$ of $a$ is given by:  
\begin{description}
 \item[ \it $a$ is 1-valent, $\Gamma_1'$ has only one true vertex]
\begin{align*}
\sigma_{\Pu^+}(a)&=  \sum_{s \in \Sw^+} \sigma_{\sw}(s) + \frac{(k-1)k}{2} (n+1) +1  \\
&\quad{} + k \cdot [n \cdot \dim W^u(a) + \dim W^u(a)] \\
&\quad{}+ [1+\dim W^u(a)]\cdot[n + |\mu(a)|], 
\end{align*}
\item[\it $a$ is 1-valent, $\Gamma_1'$ has more than one true vertex]
\begin{align*}
\sigma_{\Pu^+}(a)&= 
 \sum_{s \in \Sw^+} \sigma_{\sw}(s) + \frac{(k-1)k}{2} (n+1) +\bm_{\min}(\Gamma)+1  \\
&\quad{} + k \cdot [n \cdot \dim W^u(a) + \dim W^u(a)+1] \\
&\quad{} + \dim W^u(a)\cdot[n + |\mu(a)|+  \bm(\Gamma_1')+1]\\
&\quad{} + [n+ |\mu(a)|] \cdot[\bm(\Gamma_1') + 1]  + n +\dim \Ker_p (\Gamma'_1),
\end{align*}
\item[\it $a$  is 2-valent]
\begin{align*}
\sigma_{\Pu^+}(a)&= 
\sigma_{\pos,2}(a) + |\mu(p_1)|\cdot[\bm(\Gamma'_1) + \bm(\Gamma'_2) + |\mu(a)|] +\bm(\Gamma_1')\cdot[n+1]  \\
&\quad{} +\bm(\Gamma'_2)\cdot[n + |\mu(a)|] + e(\Gamma_1') \cdot[ \bm(\Gamma_2') +1 + e(\Gamma_2')] \\
&\quad{}   + \dim \Ker_{p_1} (\Gamma'_1) + \dim \Ker_{p_2} (\Gamma'_2)\\
&\quad{}+\bm_{\min}(\Gamma_1')+\bm_{\min}(\Gamma_2')+\bm(\Gamma_2').
\end{align*}
\end{description}

\subsection{Definition of \texorpdfstring{$\nu_{\stab}$}{v}}
With the same notation as above we are now ready to give the definition of $\nu_{\stab}$.
\begin{definition}\label{def:stab}
Let $\Gamma$ be a rigid flow tree with positive puncture $a$. 
Then the sign $\nu_{\stab}$ is given by $\nu_{\stab} = (-1)^{\mu_{\stab}}$, where
\begin{equation*}
\mu_{\stab} = \sigma_{\Pu^+}(a) + \sum_{p \in \Pu_1^-}\sigma_{\Pu_1^-}(p)  + \sum_{p\in \Pu_{2}^-}\sigma_{\Pu_{2}^-}(p)   
+\sum_{p\in \Sw^-} \sigma_{\Sw^-}(p) + 
\sum_{p\in \Y_0}\sigma_{\Y_0}(p)+ \sum_{p\in \Y_1}\sigma_{\Y_1}(p).
\end{equation*}
\end{definition}

%

\section{Examples}\label{sec:ex}
\begin{figure}[ht]
      \labellist
\small\hair 2pt
\pinlabel $a_1$ [Br] at 230 115
\pinlabel $a_2$ [Br] at 315 125
\pinlabel $a_3$ [Br] at  330 195 
\pinlabel $a_4$ [Br] at 400 190
\pinlabel $a_5$ [Br] at 420 115
\pinlabel $a_6$ [Br] at 515 219
\pinlabel $a_7$ [Br] at 515 54
\pinlabel $\partial_{x_1}$ [Br] at 190 15
\endlabellist
 \vspace{-10.8cm}
\centering
 \includegraphics[width=11cm, height=15cm]{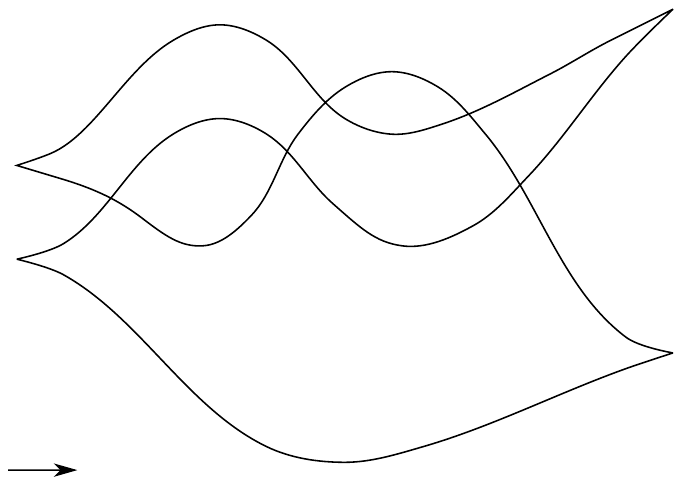} 
 \hspace{6.5cm}
 \caption{Front of $\Lambda \subset J^1(\R) = \R^3$ with orientation of $\R$ indicated.}
 \label{fig:ng}
 \end{figure}

In this section we compute the signs of the Legendrian knot with front as given in Figure \ref{fig:ng}, and we compare the result with the result we get by using the algorithm from [\cite{ngcomp}, Section 2].   

\begin{figure}[ht]
      \labellist
\small\hair 2pt
\pinlabel $e_1$ [Br] at 70 200
\pinlabel $a_3$ [Br] at 165 200
\pinlabel $e_1$ [Br] at  277 200 
\pinlabel $a_1$ [Br] at 375 200
\pinlabel $a_2$ [Br] at 470 200
\pinlabel $a_3$ [Br] at 565 200
\pinlabel $a_2$ [Br] at 540 118
\pinlabel $a_1$ [Br] at 540 88
\pinlabel $e_1$ [Br] at  540 59 
\pinlabel $\Gamma_{31}$ [Br] at 20 166
\pinlabel $\Gamma_{32}$ [Br] at 620 166 
\endlabellist
 \vspace{-7.8cm}
\centering
 \includegraphics[width=8cm, height=11cm]{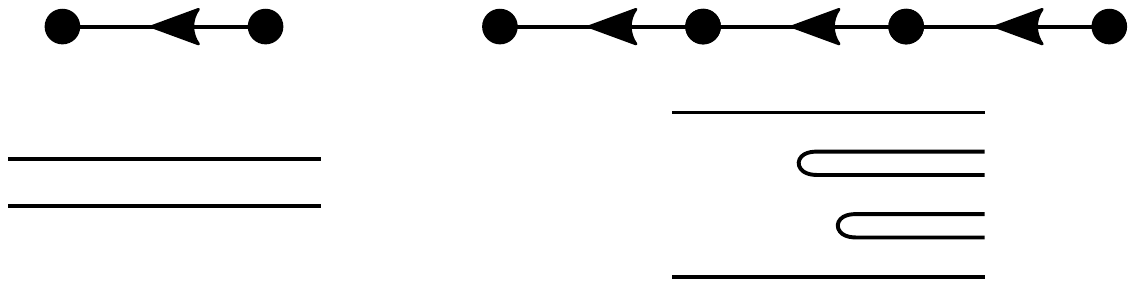} 
 \caption{Trees with $a_3$ as positive puncture. The vertices $a_1$ and $a_2$ are 2-valent negative punctures of Type 2.}
 \label{fig:t3}
 \end{figure}

\begin{figure}[ht]
      \labellist
\small\hair 2pt
\pinlabel $e_1$ [Br] at 60 100
\pinlabel $a_1$ [Br] at 150 100
\pinlabel $a_4$ [Br] at  250 100 
\pinlabel $a_5$ [Br] at 350 100
\pinlabel $a_1$ [Br] at 615 100
\pinlabel $e_1$ [Br] at 615 70
\pinlabel $a_4$ [Br] at 615 39
\pinlabel $\Gamma_{51}$ [Br] at 15 72
\endlabellist
 \vspace{-8.8cm}
\centering
 \includegraphics[width=8cm, height=11cm]{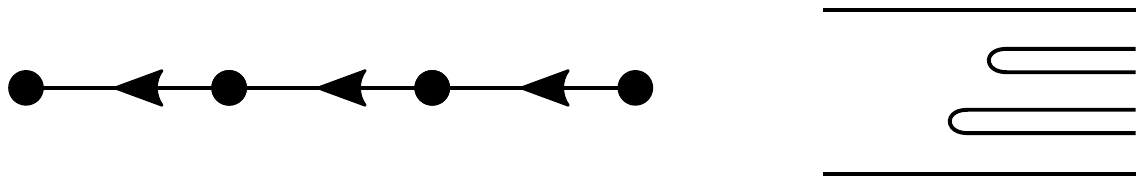} 
  \caption{The only tree with $a_5$ as positive puncture. The vertices $a_1$ and $a_4$ are 2-valent negative punctures, where $a_1$ is of Type 2 and $a_4$ is of Type 1.}
 \label{fig:t5}
 \end{figure}

\begin{figure}[ht]
      \labellist
\small\hair 2pt
\pinlabel $a_6$ [Br] at 60 195
\pinlabel $e_2$ [Br] at 169 195
\pinlabel $e_1$ [Br] at 340 195
\pinlabel $a_1$ [Br] at 435 195
\pinlabel $a_6$ [Br] at  525 195 
\pinlabel $a_1$ [Br] at 535 110
\pinlabel $e_1$ [Br] at 535 77
\pinlabel $\Gamma_{61}$ [Br] at 25 165
\pinlabel $\Gamma_{62}$ [Br] at 580 165
\endlabellist
 \vspace{-7.8cm}
\centering
 \includegraphics[width=8cm, height=11cm]{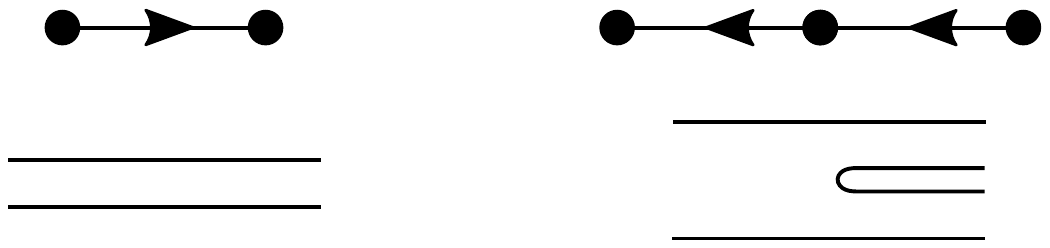} 
 \caption{Trees with $a_6$ as positive puncture. The vertex $a_1$ is a 2-valent negative puncture of Type 2.}
 \label{fig:t6}
 \end{figure}

\begin{figure}[ht]
      \labellist
\small\hair 2pt
\pinlabel $a_7$ [Br] at 250 730
\pinlabel $e_4$ [Br] at 355 730
\pinlabel $\Gamma_{71}$ [Br] at 205 700
\pinlabel $e_3$ [Br] at 50 580
\pinlabel $a_1$ [Br] at  145 580
\pinlabel $a_7$ [Br] at  245 580
\pinlabel $e_3$ [Br] at 355 580
\pinlabel $a_5$ [Br] at  445 580
\pinlabel $a_7$ [Br] at  545 580
\pinlabel $e_3$ [Br] at 245 498
\pinlabel $a_1$ [Br] at  245 460
\pinlabel $e_3$ [Br] at 545 498
\pinlabel $a_5$ [Br] at  545 460
\pinlabel $\Gamma_{72}$ [Br] at 5 550
\pinlabel $\Gamma_{73}$ [Br] at 590 550
\pinlabel $e_3$ [Br] at 115 395
\pinlabel $a_1$ [Br] at  208 395 
\pinlabel $a_3$ [Br] at  305 395 
\pinlabel $a_4$ [Br] at  405 395 
\pinlabel $a_7$ [Br] at  508 395 
\pinlabel $e_3$ [Br] at 405 315
\pinlabel $a_1$ [Br] at  405 288
\pinlabel $a_3$ [Br] at  405 266 
\pinlabel $a_4$ [Br] at  405 240
\pinlabel $\Gamma_{74}$ [Br] at 60 365 
\pinlabel $e_3$ [Br] at 110 193
\pinlabel $a_1$ [Br] at  205 193 
\pinlabel $a_2$ [Br] at  305 193 
\pinlabel $a_5$ [Br] at  399 193 
\pinlabel $a_7$ [Br] at  495 193 
\pinlabel $e_3$ [Br] at 405 113
\pinlabel $a_1$ [Br] at  405 89
\pinlabel $a_2$ [Br] at  405 63 
\pinlabel $a_5$ [Br] at  405 35
\pinlabel $\Gamma_{75}$ [Br] at 60 165 
\endlabellist
 \vspace{-.8cm}
\centering
 \includegraphics[width=8cm, height=11cm]{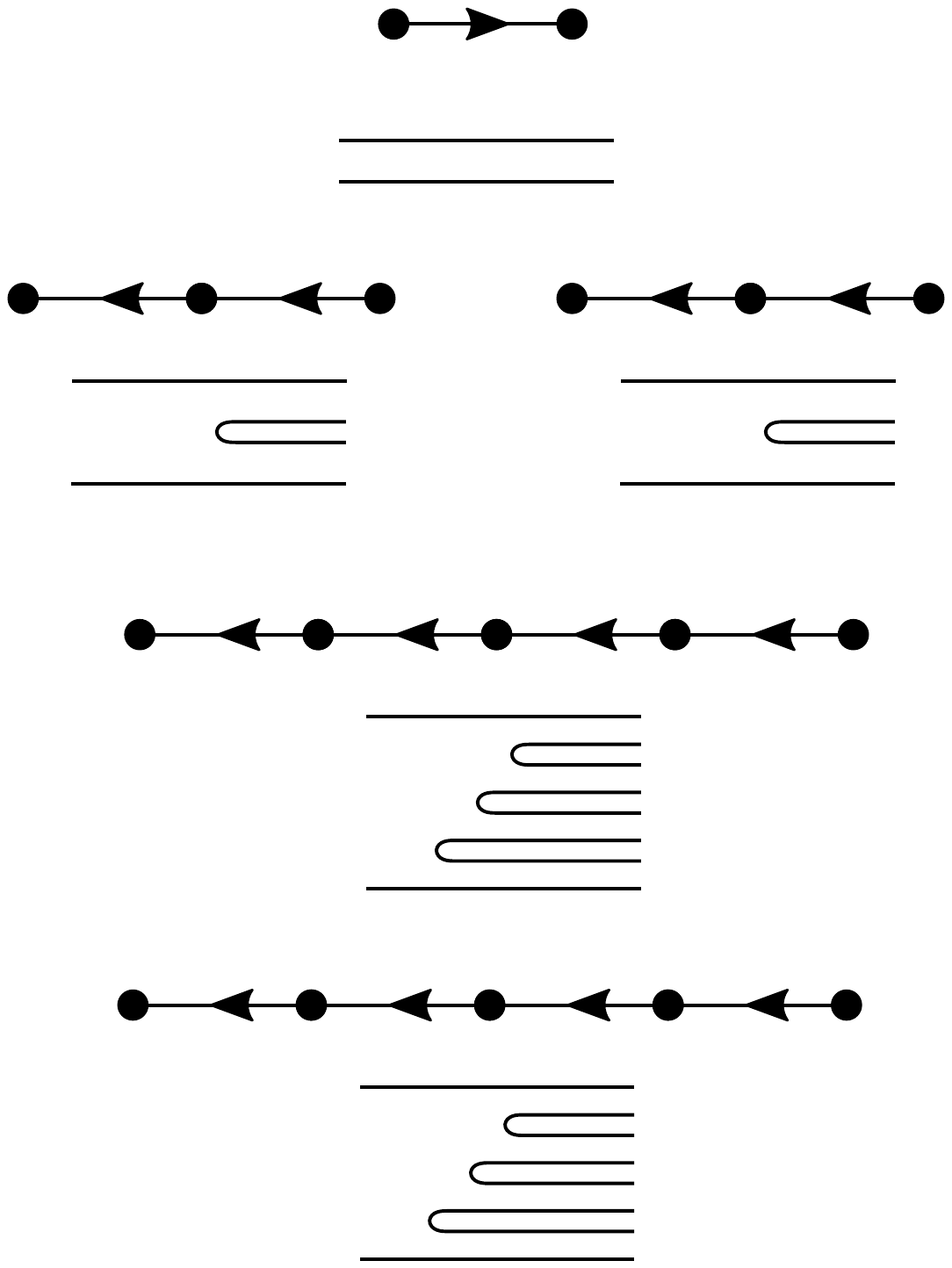} 
 \caption{Trees with $a_7$ as positive puncture. The vertices $a_1, a_2, a_3, a_4, a_5$ are 2-valent negative punctures of Type 1.}
 \label{fig:t7}
 \end{figure}

We get the following differential, modulo 2:
\begin{align*}
 &\partial a_1 = \partial a_2 =\partial a_4 = 0 \\
 & \partial a_3 = 1 + a_1a_2 \\
 & \partial a_5 = a_4a_1 \\
 & \partial a_6 =1 + a_1 \\
 & \partial a_7 =1 + a_1 + a_5 +  a_5a_2a_1 + a_4a_3a_1.
\end{align*}
Indeed, it is easy to see that we do not have any flow trees with $a_1, a_2$ or $a_4$ as positive puncture. The flow trees for the other punctures are pictured in Figures \ref{fig:t3} -- \ref{fig:t7}, together with their standard domains. Here $e_1$ corresponds to the end in the upper left corner of the picture of the knot, $e_2$ is the end in the upper right corner, $e_3$ is the end in the lower left corner and $e_4$ is the end in the lower right corner. It remains to calculate the signs of the trees.

\begin{figure}[ht]
      \labellist
\footnotesize\hair 1.5pt
\pinlabel $r_d(-\pi)$ [Br] at 70 420
\pinlabel $(\Lambda'(\theta),\partial_{x_2})$ [Br] at 10 330
\pinlabel $(\partial_{x_1},-\partial_{x_2})_d$ [Br] at 180 400
\pinlabel $(-\partial_{x_1},\partial_{x_2})_s$ [Br] at 180 350
\pinlabel  $r_s(-\pi)$ [Br] at 245 410
\pinlabel $r_d(-\pi)$ [Br] at 565 280
\pinlabel $r_s(\pi)$ [Br] at 512 302
\pinlabel $r_s(-\pi)$ [Br] at 508 400
\pinlabel $r_s(\pi)$ [Br] at 408 390
\pinlabel $(\Lambda'(\theta),\partial_{x_2})$ [Br] at 680 330
\pinlabel $(\partial_{x_1},\partial_{x_2})$ [Br] at 355 410
\pinlabel $r_s(-\pi)$ [Br] at 368 333
\pinlabel $r_s(\pi)$ [Br] at 328 313
\pinlabel $r_d(-\pi)$ [Br] at 70 -10
\pinlabel $(\Lambda'(\theta),\partial_{x_2})$ [Br] at 15 64
\pinlabel $(\partial_{x_1},-\partial_{x_2})_d$ [Br] at 300 42
\pinlabel $(-\partial_{x_1},\partial_{x_2})_s$ [Br] at 300 -10
\pinlabel $r_s(-\pi)$ [Br] at 443 -10
\pinlabel $(\Lambda'(\theta),\partial_{x_2})$ [Br] at 670 50
\pinlabel $r_d(-\pi)$ [Br] at 610 110
\pinlabel $r_s(\pi)$ [Br] at 510 110
\pinlabel $(\partial_{x_1},\partial_{x_2})$ [Br] at 380 120
\endlabellist
 \vspace{-4.8cm}
 \hspace{-.8cm}
\centering
 \includegraphics[width=12cm, height=14cm]{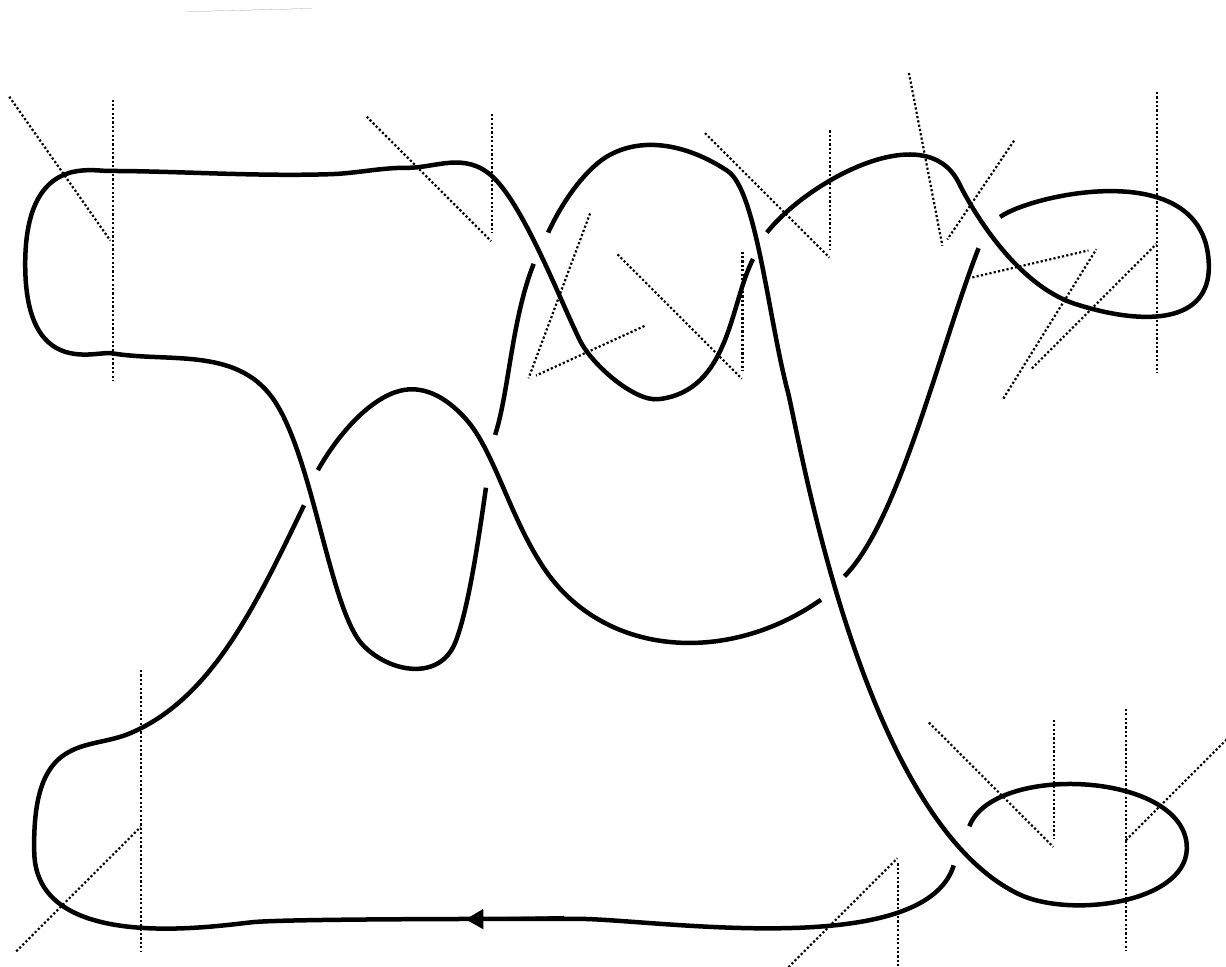} 
 \caption{Lagrangian projection of $\Lambda$ together with trivializations $\tau_d$ and $\tau_s$. The orientation of $\Lambda$ is indicated by an arrow and $\partial_{x_1}$ points to the right in the picture. The subscript $d$ means  the default trivialization and the subscript $s$ means the trivialization induced by the spin structure, and $r(\pm \pi)$ is a $\pm\pi$-rotation in the $\Lambda'(\theta), \partial_{x_2}$-plane. In the case the subscript is missing the two trivializations agree. }
 \label{fig:nglag}
 \end{figure}

First we fix the initial orientation choices. We choose the orientation of $\Lambda$ as in Figure \ref{fig:nglag} and of $M = \R$ as in Figure \ref{fig:ng}.  We pick the Lie group spin structure on $\Lambda$ and we fix an orientation of $\C$ as in [\cite{orienttrees}, Section 9.1.B]. Finally we choose the orientation of all the relevant unstable manifolds to be given by $\partial_{x_1}$, except for the ones associated to $a_6$ and $a_7$ where we choose the orientation to be given by $-\partial_{x_1}$. We also note that 
\begin{equation*}
|\mu(a_1)|=|\mu(a_2)|=|\mu(a_5)|=0, \quad |\mu(a_3)|=|\mu(a_4)|=|\mu(a_6)|=|\mu(a_7)|=1.
\end{equation*}  


Next consider the trivialization of the stabilized tangent bundle $T \tilde \Lambda$ along parts corresponding to vertex regions of the trees of $\Lambda$ as described in Section  \ref{sec:trivtriv}. Clearly these trivializations patch together to a trivialization $\tau_d$ of $T \tilde \Lambda$, using the identity to glue together the trivializations from different vertex regions. See Figure \ref{fig:nglag}, where $\Lambda(\theta)$, $\theta \in S^1$, is a parameterization of $\Lambda$ so that $\Lambda'(\theta)$ gives the positive orientation of $T_{\Lambda(\theta)}\Lambda$, and $r
_d(\pm \pi)$ gives a positive/negative $\pi$-rotation in the $(\Lambda'(\theta),\partial_{x_2})$-plane as described in Section \ref{sec:trivtriv}. The subscript $d$ indicates that we are considering the default trivialization. We also have a trivialization $\tau_s$ induced by the spin structure. To describe this one, start at the arrow in Figure \ref{fig:nglag}, which indicates the chosen orientation of $\Lambda$. Here the spin trivialization is given by $(\Lambda'(\theta),\partial_{x_2}) = (-\partial_{x_1},\partial_{x_2})$. Follow the knot in the direction of its orientation, and perform $\pm\pi$-rotations in the $(\Lambda'(\theta),\partial_{x_2})$-plane in the regions marked by $r_s(\pm\pi)$. These rotations are homotopies of the trivialization $(\Lambda'(\theta),\partial_{x_2})$ that we introduce so that the $\tau_d$- and $\tau_s$-trivialization coincide at Reeb chord end points. Away from these rotations we let the trivialization be given by $\pm(\Lambda'(\theta),\partial_{x_2})$ (note that a $\pi$-rotation changes the trivialization from  $\pm(\Lambda'(\theta),\partial_{x_2})$ to $\mp(\Lambda'(\theta),\partial_{x_2})$, and similarly for a $-\pi$-rotation). 

\begin{remark}
If we make another choice of homotopy some signs in the computation below may change, but this can be compared with making different choices of the orientations of the stable and unstable manifolds.
\end{remark}

\begin{lemma}\label{lemma:spintriv}
 We have $\nu_{\triv}(\Gamma)=-1$ for any flow tree $\Gamma$ of $\Lambda$ containing $e_2$ or $e_4$, and $\nu_{\triv}(\Gamma)=1$ otherwise.
\end{lemma}
\begin{proof}
Consider the trees $\Gamma_{61}$ and $\Gamma_{71}$. The cotangent lift of these trees consist of one arc each, and we see that if we compare $\tau_s$ and $\tau_d$ along each of those arcs we get a homotopically nontrivial loop in $SO(2)$, in both cases. Hence $\nu_{\triv}(\Gamma_{61})=\nu_{\triv}(\Gamma_{71})=-1$.

For the other trees we see that along each of the connected components of the cotangent lifts of the trees the trivializations $\tau_d$ and $\tau_s$ are equal, up to homotopy.  That is, either we have cancelling pairs $r_s(\pm\pi), r_s(\mp \pi)$ or we have $r_s(-\pi)$ and $r_d(-\pi)$ along the same arc. This implies that $\nu_{\triv}(\Gamma)=1$ for these trees.
\end{proof}

 Now we compute the remaining signs for the trees $\Gamma_{31}$ and $\Gamma_{32}$ in Figure \ref{fig:t3}. We start with $\Gamma_{13}$, which we divide into two sub flow trees $\Gamma_e$ and $\Gamma_3$ at a point $q \in \Gamma_{13}$ so that $\Gamma_e$ is the sub flow tree with the end as true vertex and $q$ as positive special puncture, and $\Gamma_3$ has $a_3$ as true puncture and $q$ as a negative special puncture. 
 
 Clearly we have 
 \begin{equation*}
  \nu_{\en}(\Gamma_{13}) = 1.
 \end{equation*}
To compute $\nu_{\inter}(\Gamma_{13})$ we see that with our conventions from Section \ref{sec:flowor} we have 
\begin{align*}
 & \Or(\F_q(\Gamma_3)) = \partial_{x_1} \\
 & \Or(\F_q(\Gamma_e)) = \partial_{x_1}.
\end{align*}
Since $\Gamma_{13}$ only have two vertices we see that $\nu_{\inter}(\Gamma_{13})$ should be computed using Case 1a) in Section \ref{sec:finglu}. The flow direction of $\Gamma_{31}$ is given by $- \partial_{x_1}$, and 
\begin{equation*}
 \Or_{\inter}(\F_q(\Gamma_e), \F_q(\Gamma_3)) = \partial_{x_1}.
\end{equation*}
Hence 
\begin{equation*}
 \nu_{\inter}(\Gamma_{13}) = -1.
\end{equation*}
It remains to compute $\nu_{\stab}(\Gamma_{31})=(-1)^{\sigma_{\Pu^+}(a_3)}$ where in this case 
\begin{align*}
 \sigma_{\Pu^+}(a_3) &\equiv 0 + 0 + 1 + 0 + [1 + 1] \cdot [1 + |\mu(a_3)|]\\
 &\equiv  1, \pmod{2},
\end{align*}
since $|\mu(a_3)| =1 $
 
Thus we get 
\begin{align*}
 \nu(\Gamma_{13}) &= 1 \cdot (-1) \cdot 1 \cdot (-1)^{1} \\
 &=  1.
\end{align*}

To compute $\nu(\Gamma_{32})$, we divide the tree into four sub flow trees $\Gamma_e, \Gamma_1,$ $\Gamma_2, \Gamma_3$, where $\Gamma_e$ contains the end as only true puncture, and $\Gamma_i$ contains the puncture $a_i$ as only true puncture, $i=1,2,3$. Let $q_i$ denote the negative special puncture of $\Gamma_i$, $i=1,2,3$. Let $\Gamma_{e1}$ be the sub flow tree with $e,a_1$ as true punctures and $q_2$ as special positive puncture, and let $\Gamma_{e12}$ be the sub flow tree we get from gluing $\Gamma_{e1}$ to $\Gamma_2$. Then 
\begin{align*}
 & \Or(\F_{q_3}(\Gamma_3)) = \partial_{x_1} \\
 & \Or(\F_{q_1}(\Gamma_e)) = \partial_{x_1},
\end{align*}
and 
\begin{equation*}
 \Or(\I_{a_1}(\Gamma_{e1})) = \Or_{\inter}(\partial_{x_1}, {\bf{1}}) = (-1)^{0 + 1\cdot 1}{\bf{1}} = (-1)\cdot {\bf{1}},
\end{equation*}
where ${\bf 1}$ denotes the positively oriented 0-dimensional vector space. Hence 
\begin{equation*}
 \Or(\F_{q_2}(\Gamma_{1e})) = - \partial_{x_1},
\end{equation*}
and so 
\begin{equation*}
 \Or(\I_{a_2}(\Gamma_{e12})) = \Or_{\inter}(-\partial_{x_1}, {\bf{1}}) = {\bf{1}}.
\end{equation*}
From Section \ref{sec:finglu}, Case 1b), it follows that $\nu_{\inter}(\Gamma_{32}) $ satisfies 
\begin{equation*}
\nu_{\inter}(\Gamma_{32}){\bf 1} \wedge \partial_{x_1} = \partial_{x_1},
\end{equation*}
and hence 
\begin{equation*}
 \nu_{\inter}(\Gamma_{32}) = 1.
\end{equation*}

Next we compute $\nu_{\en}(\Gamma_{32})$, which is given by the equation 
\begin{equation*}
 \nu_{\en}(\Gamma_{32}) e \wedge a_1 \wedge a_2 = a_1 \wedge a_2 \wedge e,
\end{equation*}
and so 
\begin{equation*}
 \nu_{\en}(\Gamma_{32}) = 1.
\end{equation*} 
 
Finally we compute $\nu_{\stab}(\Gamma_{32}) = (-1)^{\sigma_{\Pu^+}(a_3) + \sigma_{\Pu_{2}^-}(a_1) + \sigma_{\Pu_{2}^-}(a_2)}$. We have
\begin{align*}
 \sigma_{\Pu^+}(a_3) &=  0+0+2 + 1 + 0  \\
 &\quad{}  + 1\cdot[1+|\mu(a_3)| + 2 +1] + 
  [1 + |\mu(a_3)|]\cdot[2+1] +1 +0 \\
 &\equiv 1, \pmod{2}, \\
 \sigma_{\Pu_{2}^-}(a_1) &= \sigma_{\negg,2}((k_2,k'_2),|\mu(a_1)|, 1, 2) + 0 + 0 + 1 +1 + 0 \\
 & \equiv \sigma_{\negg,2}((1,0),0, 1, 2), \pmod{2}, \\
 \sigma_{\Pu_{2}^-}(a_2) &= \sigma_{\negg,2}((k_3,k'_3),|\mu(a_2)|, |\mu(q_2)|, 2) + 1 + 1 + 0 +1 + 1\cdot[1 + |\mu(a_2)| +1] \\
 & \equiv \sigma_{\negg,2}((1,0),0, 1, 2) + 1, \pmod{2},
\end{align*}
where $(k_i,k'_i)$ is the sheet signs of $q_i$, $i=2,3$. 

Hence we get 
\begin{align*}
 \nu_{\stab}(\Gamma_{32}) = (-1)^{1+\sigma_{\negg,2}((1,0),0, 1, 2) + \sigma_{\negg,2}((1,0),0, 1, 2) + 1},
\end{align*}
and thus 
\begin{equation*}
 \nu(\Gamma_{32}) =1. 
\end{equation*} 

In a similar way one computes the signs of the other trees. We have summarized the result in Table \ref{table:treesigns}, where we use the notation $\nu_{\triv} = (-1)^{\mu_{\triv}}$.

  \begin{table}[hhh]
 \begin{center}
   \caption{{\bf Signs for the trees of the knot $\Lambda$.
   }}\label{table:treesigns}
 \tabulinesep=1.15mm
  \begin{tabu}{|>{\arraybackslash}m{.3in} | >{\arraybackslash}m{.3in} | >{\arraybackslash}m{.3in} | >{\arraybackslash}m{.3in} | >{\arraybackslash}m{1.5in} | >{\arraybackslash}m{1.5in} |}
  \hline  
  Tree & $\mu_{\triv}$ & $\mu_{\inter}$ & $\mu_{\en}$ & $\mu_{\stab} $ & $\mu_{\triv}+\mu_{\inter}+\mu_{\en}+\mu_{\stab}$\\
   \hline 
    $\Gamma_{31}$ & $0$ & $1$ & $0$ & $1$ & 0 \\
    \hline
    $\Gamma_{32}$ & $0$ & $0$ & $0$ & $0$ & $0$\\
    \hline
    $\Gamma_{51}$ & $0$ & $0$ & $1$ & $\sigma_{\negg,2}((1,0),0, 1, 2) + \sigma_{\negg,2}((0,0),1, 1, 1)$ & $\sigma_{\negg,2}((1,0),0, 1, 2) + \sigma_{\negg,2}((0,0),1, 1, 1)+ 1$\\
    \hline
    $\Gamma_{61}$ & $1$ & $0$ & $0$ & $1$ & $0$\\
    \hline
    $\Gamma_{62}$ & $0$ & $0$ & $1$ & $\sigma_{\negg,2}((1,0),0, 1, 2) +1$ & $\sigma_{\negg,2}((1,0),0, 1, 2)$\\
    \hline
     $\Gamma_{71}$ & $1$ & $0$ & $0$ & $ 1$ & $0$\\
    \hline
    $\Gamma_{72}$ & $0$ & $0$ & $0$ & $\sigma_{\negg,2}((0,1),0, 1, 1) )$ & $\sigma_{\negg,2}((0,1),0, 1, 1)$\\ 
    \hline
    $\Gamma_{73}$ & $0$ & $0$ & $0$ & $\sigma_{\negg,2}((0,1),0, 1, 1) $ & $\sigma_{\negg,2}((0,1),0, 1, 1) $\\
    \hline
     $\Gamma_{74}$ & $0$ & $0$ & $0$ & $\sigma_{\negg,2}((0,1),0, 1, 1) + \sigma_{\negg,2}((0,1),1, 0, 1)+ \sigma_{\negg,2}((1,1),1, 1, 1)+ 1$ & $\sigma_{\negg,2}((0,1),0, 1, 1) + \sigma_{\negg,2}((0,1),1, 0, 1)+ \sigma_{\negg,2}((1,1),1, 1, 1)+1$\\
    \hline
     $\Gamma_{75}$ & $0$ & $0$ & $0$ & $\sigma_{\negg,2}((0,1),0, 1, 1) $ & $\sigma_{\negg,2}((0,1),0, 1, 1) $\\
    \hline
\end{tabu}
 \end{center}
\end{table}

Thus we get 
\begin{align*}
 \partial a_1 &= \partial a_2 =\partial a_4 = 0 \\
  \partial a_3 &=1 + a_1a_2 \\
  \partial a_5 &= (-1)^{\sigma_{\negg,2}((1,0),0, 1, 2) + \sigma_{\negg,2}((0,0),1, 1, 1)+ 1}a_4a_1 \\
  \partial a_6 &=1 + (-1)^{\sigma_{\negg,2}((1,0),0, 1, 2)} a_1 \\
  \partial a_7 &=1 + (-1)^{\sigma_{\negg,2}((0,1),0, 1, 1)}(a_1 + a_5) \\
 &\quad{} +  (-1)^{\sigma_{\negg,2}((0,1),0, 1, 1) + \sigma_{\negg,2}((0,1),1, 0, 1)+ \sigma_{\negg,2}((1,1),1, 1, 1)+1}a_4a_3a_1 \\
 &\quad{} + (-1)^{\sigma_{\negg,2}((0,1),0, 1, 1)} a_5a_2a_1.
\end{align*}

If we compute the differential using the algorithm in  \cite{ngcomp} we get 
\begin{align*}
 &\partial a_1 = \partial a_2 =\partial a_4 = 0 \\
 & \partial a_3 = 1 + a_1a_2 \\
 & \partial a_5 = - a_4a_1 \\
 & \partial a_6 =1 + a_1 \\
 & \partial a_7 =1 - a_1 - a_5 -  a_5a_2a_1 - a_4a_3a_1.
\end{align*}

 Thus, if we have 
\begin{align*}
 &\sigma_{\negg,2}((1,0),0, 1, 2) \equiv 0 \\
  &\sigma_{\negg,2}((0,0),1, 1, 1)  \equiv 0 \\
 &\sigma_{\negg,2}((0,1),0, 1, 1)\equiv 1 \\
 &\sigma_{\negg,2}((0,1),1, 0, 1)+ \sigma_{\negg,2}((1,1),1, 1, 1) \equiv 1, \\
\end{align*}
modulo 2, we get the same signs as from \cite{ngcomp}.

\begin{remark}
In this case, we do not know if the orientation of $\C$ agrees with the one chosen in [\cite{orienttrees}, Section 9], and hence we cannot compare with our formulas for $\sigma_{\negg,2}$.
\end{remark}

\newpage
\bibliographystyle{halpha}
\bibliography{c_main_n} 
%
%

\end{document}